\documentclass[preprint]{elsarticle}
\usepackage[utf8]{inputenc}
\usepackage[T1]{fontenc}
\usepackage[french,english]{babel}
\usepackage{amssymb}
\usepackage{amsmath}
\usepackage{graphicx}
\usepackage{algorithm}
\usepackage[noend]{algorithmic}
\usepackage{vmargin}
\usepackage{url}
\usepackage{tikz}
\usepackage{natbib}
\usepackage{xcolor,colortbl}
\usepackage{rotating}

\usetikzlibrary{calc}
\usetikzlibrary{plotmarks}
\usetikzlibrary{shapes}
\usetikzlibrary{decorations}
\usetikzlibrary{decorations.pathmorphing}
\usetikzlibrary{decorations.pathreplacing}
\usetikzlibrary{decorations.shapes}

\urldef{\mailENS}\path|nbourgeo@phare.normalesup.org|
\urldef{\mailSAMM}\path|marie.cottrell,patrick.letremy@univ-paris1.fr|
\urldef{\mailLAMOP}\path|benjamin.deruelle,stephane.lamasse@univ-paris1.fr|

\newcommand{\algo}[1]{\textsc{Algorithm~#1}}

\definecolor{gris}{gray}{0.6}
\definecolor{grisclair}{gray}{0.9}
\definecolor{gray0.9}{gray}{.95}
\definecolor{gray0.8}{gray}{.85}
\definecolor{gray0.7}{gray}{.75}
\definecolor{gray0.6}{gray}{.7}
\definecolor{gray0.55}{gray}{.65}
\definecolor{gray0.5}{gray}{.6}
\definecolor{gray0.45}{gray}{.575}
\definecolor{gray0.4}{gray}{.55}
\definecolor{gray0.35}{gray}{.525}
\definecolor{gray0.3}{gray}{.5}

\makeatletter

\let\c@table\c@figure
\makeatother 

\journal{Neurocomputing}

\begin{document}
\begin{frontmatter}

\title{How to improve robustness in Kohonen maps and visualization in Factorial Analysis: application to text mining}

\author[samm]{Nicolas Bourgeois}
\author[samm]{Marie Cottrell}
\author[lamop]{Benjamin D\'eruelle}
\author[lamop]{St\'ephane Lamass\'e}
\author[samm]{Patrick Letr\'emy}

\address[samm]{SAMM - Universit\'e Paris 1 Panth\'eon-Sorbonne 90, rue de Tolbiac, 75013 Paris, France \mailENS,\mailSAMM}
\address[lamop]{PIREH-LAMOP - Universit\'e Paris 1 Panth\'eon-Sorbonne 1, rue Victor Cousin, Paris, France \mailLAMOP}

%\maketitle

\begin{abstract}
This article is an extended version of a paper presented in the WSOM'2012 conference \cite{wsom12}. We display a combination of factorial projections, SOM algorithm and graph techniques applied to a text mining problem. The corpus contains 8 medieval manuscripts which were used to teach arithmetic techniques to merchants. 

Among the techniques for Data Analysis, those used for Lexicometry (such as Factorial Analysis) highlight the discrepancies between manuscripts. The reason for this is that they focus on the deviation from the independence between words and manuscripts. Still, we also want to discover and characterize the common vocabulary among the whole corpus.

Using the properties of stochastic Kohonen maps, which define neighborhood between inputs in a non-deterministic way, we highlight the words which seem to play a special role in the vocabulary. We call them fickle and use them to improve both Kohonen map robustness and significance of FCA visualization. Finally we use graph algorithmic to exploit this fickleness for classification of words.
\end{abstract}

\end{frontmatter}

\section*{Introduction}
\setcounter{footnote}{0}

\subsection*{Historical Context}

One approach to understand the evolution of science is the study of the evolution of the language used in a given field. That is why we would like to pay attention to the vernacular texts dealing with practical arithmetic and written for the instruction of merchants. Such texts are known since the XIII\textsuperscript{th} century, and from that century onwards,  the vernacular language appears more and more as the medium of practical mathematics.\\

Treaties on arithmetical education were therefore mostly thought and written in local languages, (they  were written not only in French but also in Italian, Spanish, English and German). In this process, the XV\textsuperscript{th} century appears as a time of exceptional importance because we can study the inheritance of two hundred years of practice. For the authors of these texts, the purpose was not only to teach merchants, but also to develop knowledge in vernacular language. Their books were circulated far beyond the shopkeepers' world, to the humanists' circles for example.\\

\subsection*{An objective of historical research: the study of specialized languages}

The work previously done by historians~\cite{lamasse2012} consisted in the elaboration of a dictionary of the lexical forms found in all the treaties, in order to identify the different features of the mathematical vernacular language at that time. This being done, we have worked on the contexts of some especially important words in order to understand the lexicon in all its complexity. In other words, we would like to determine the common language that forms the specialized language beyond the specificity of each text.

\subsection*{Outline of this work}

%This work meets a double purpose: on one hand we study a specific, well-defined corpus, which is based on a set of scientific texts from the XV\up{th} century, on the other hand, we introduce a new tool, which can be used for any problem of lexicometry - or any problem involving contingency matrices.
Among the techniques for Data Analysis, those used for Lexicometry (such as Factorial Analysis) highlight the discrepancies between manuscripts. The reason for this is that they focus on the deviation from the independence between words and manuscripts. Still, we also want to discover and characterize the common vocabulary among the whole corpus. That is why we introduce a new tool, which combine the properties of Factorial Correspondence Analysis and Stochastic Self-Organizing Maps. That leads to the definition of fickle pairs and fickle words. Fickle words can be seen as this common vocabulary we are looking for, and prove themselves to be a good basis for a new visualization with the help of graph theory.

In part~\ref{corpus}, we first focus on the definition of the corpus: the texts, the pre-processing, and the protocol which is traditionally used in Humanities and Social Sciences to handle such data. Then, (part~\ref{method}), we design the tools : 'fickle pairs' and 'fickle words', robust Kohonen Maps, improved FCA, graphs of relations between words based on fickleness. We explain the algorithms involved and display the results on the corpus. Finally (part~\ref{analysis}), we give a brief analysis and comments on these results.

\section{Text, Corpus and protocol}\label{corpus}

In order to delimit a coherent corpus among the whole European production of practical calculation education books, we have chosen to pay attention to those treaties which are sometimes qualified as commercial (\textit{marchand} in French) which have been written in French between 1415 and about 1500. Note that this corpus has already been studied by~\cite{beaujouanplace1988}, \cite{spiesserarithmetique2003} and~\cite{benoitrecherches1985}. In this way, our corpus follows the rules of the discourse analysis: homogeneity, contrastiveness and diachronism. For further explanation about texts, methodology and purpose of the analysis see~\cite{lamasse2012}, for further explanation about the corpus~\cite{prost1988}, for wider explication about analysis see~\cite{Mayaffre2005}, \cite{Rastier2011}.\\

It contains eight treaties on the same topic, written in the same language and by different XV\textsuperscript{th} century authors. The following Table~\ref{tab:chronomss} describes some elements of the lexicometric characteristics of the corpus and shows how non balanced it is.

\hyphenpenalty=100000
\begin{table}
\begin{center}
\begin{small}
\setlength{\tabcolsep}{4pt}
\renewcommand{\arraystretch}{1.2}
\begin{tabular}{|p{4cm}|p{1.5cm}|p{1.5cm}|p{1.8cm}|p{1.4cm}|p{1.3cm}|}
\hline
\textbf{Manuscripts} and \textbf{Title} &   \textbf{Date} & \textbf{Author}  & \textbf{Number of} \textbf{occurrences} & \textbf{Number of} \textbf{words} & \textbf{Hapax} \\\hline\hline
Bib. nat. Fr. 1339 & ca.1460 & anonyme & 32077 & 2335 & 1229 \\
Bib. nat. Fr. 2050 & ca.1460 & anonyme & 39204 & 1391 & 544 \\
Cesena Bib. Mal. S-XXVI-6, \textit{Traict\'e de la praticque}  & 1471? & Mathieu {Pr\'ehoude}? & 70023 & 1540 & 635 \\
Bibl. nat. Fr. 1346, Commercial appendix~of \textit{Triparty en la science des nombres} & 1484		& Nicolas Chuquet & 60814 & 2256 & 948 \\ 
M\'ed. Nantes 456   & ca.1480-90		& anonyme & 50649 & 2252 & 998\\
Bib. nat. Arsenal 2904, \textit{Kadran aux marchans} & 1485 & Jean Certain & 33238 &1680 &714 \\ 
Bib. St. Genv. 3143 & 1471 & Jean Adam & 16986 & 1686 & 895\\
Bib. nat. Fr. Nv. Acq. 10259 & ca.1500 & anonyme & 25407 & 1597  &730 \\
\hline
\end{tabular}
\end{small}
\end{center}
\caption{Corpus of texts and main lexicometric features. The number of occurences is the total number of words including repetitions, the number of words is the number of distinct words, Hapax are words appearing once in a text.}
\label{tab:chronomss}
\end{table}
\hyphenpenalty=1000
\subsection{Humanities and Social Sciences traditional protocol}

Traditionally on this kind of textual data, researchers in Humanities and Social Sciences work on statistical specificity and contextual concordances, since they allow an easy discovery of the major lexical splits within the texts of the corpus, while remaining close to the meanings of the different forms.\\

Then, the factorial and clustering methods, combined with co-occurrences analysis (see \cite{martinezcontribution2003}) help us to cluster the texts without breaking the links with semantic analysis.\\

However, such a method of data processing requires a preliminary treatment of the corpus, the lemmatization~\cite{brunet2000}. It consists in gathering the different inflected forms of a given word as a single item. It allows us to work at many different levels of meaning, depending upon the granularity adopted: forms, lemma, syntax. 

We can justify this methodological choice here by its effect on the dispersion of the various forms which can be linked to the same lemma, a high degree of dispersion making the comparison between texts more difficult. It must also be remembered that in the case of medieval texts, this dispersion is increased by the lack of orthographic norms. In our case, this process has an important quantitative consequence on the number of forms in the corpus, which declines from 13516 forms to 9463, a reduction of some 30\%.\\

%The factorial analysis allowed us to establish a typology of the complete parts of the corpus, based upon all the forms. However, it can be useful to improve this global analysis with a probabilistic calculation for each component of the corpus, by using the table of the under-frequencies \cite{lebartstatistique1994}. It makes it possible to compare the parts of the corpus with each other, taking into account the occurrences of the words and their statistical specificities.\\

This process has been achieved with a particular attention to the meaning of each word in order to suppress ambiguities: a good example is the French word \textit{pouvoir} which can be a verb translated by "can" or "may", and which is also a substantive translated by "power".

Finally, to realize a clustering of the manuscripts, we have only kept the 219 words with highest frequencies. The set of words selected that way for text classification relates to mathematical aspects, such as operations, numbers and their manipulations, as well as to didactic aspects. Their higher frequencies reflect the fact that they are the language of the mathematics as they appear to be practiced in these particular texts. 

Thus, in what follows, the data are displayed in a contingency table $T$ with $I=219$ rows (the words) and $J=8$ columns (the manuscripts) so that the entry $t_{i,j}$ is the number of occurrences of word $i$ in manuscript $j$.

\subsection{Use of Factorial Correspondence Analysis (FCA)}\label{fcadef}

Factorial Correspondence Analysis is one of the factorial methods which consist in applying an orthogonal transformation to the data, to supply the user with simplified representation of high-dimensional data, as defined in 
\cite {Morrison1967}. The most popular of these factorial methods is the Principal Component Analysis, which deals with real-valued variables and supplies for example the best two-dimensional representation of high-dimensional dataset, by retaining the first two eigenvectors of the covariance matrix. 

Factorial Correspondence Analysis (see~\cite {Benzecri1992} or \cite {Lebart1984}) is a variant of Principal Component Analysis, designed to deal with categorical variables. Let us consider two categorical variables with respectively $I$ and $J$ items  and the associated contingency table $T$ where entry $t_{i,j}$ is the number of co-occurrences of item $i$ for the row variable and item $j$ for column variable. The rows and the columns are scaled to sum to 1 and normalized in order to be treated simultaneously, by defining

\begin{equation}
t_{i,j}^{norm}=\frac{t_{i,j}}{\sqrt{\sum_i t_{i,j} \sum_j t_{i,j}}}.\label{normalization}
\end{equation}

To achieve the FCA, two Principal Component Analysis are made over the normalized table $T^{norm}$ (providing a representation of the rows) and its transposed table (providing a representation of the columns). The main property of FCA is that  both representations can be superposed, since their principal axes are strongly correlated. The proximity between items is significant, regardless they stand for  row items or column items, except in the center of the map.

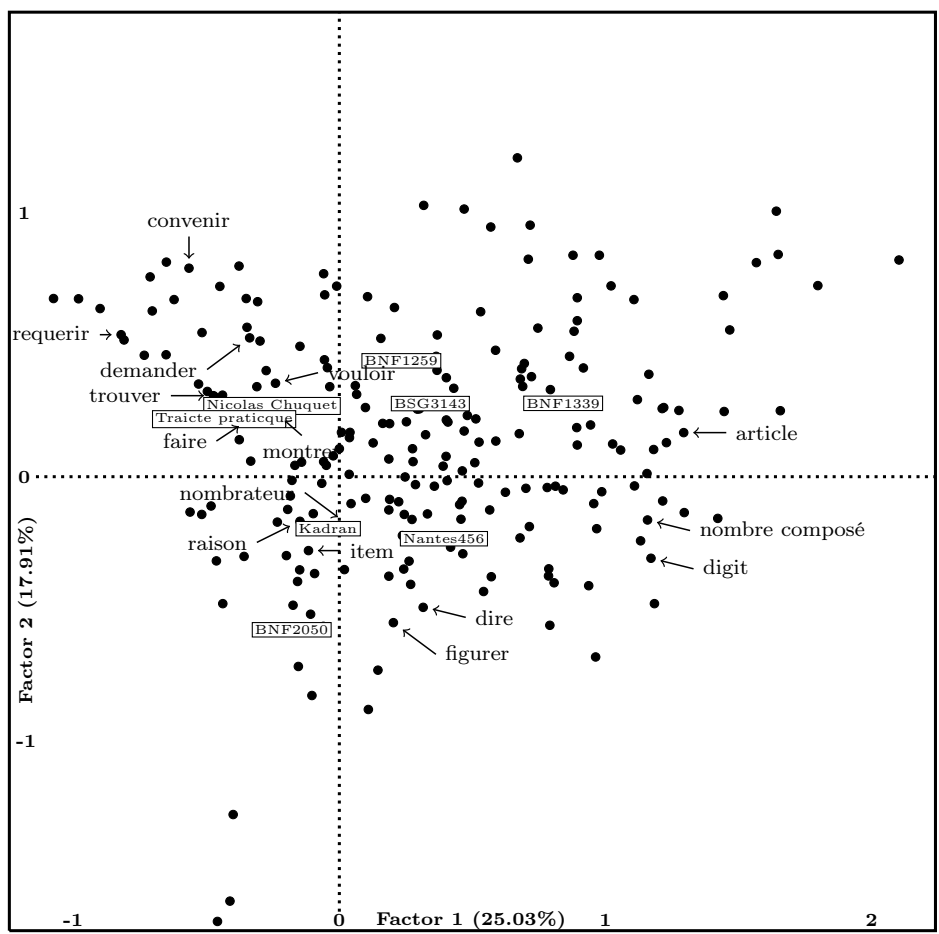
\begin{figure}[h!]
\centering

\begin{tikzpicture}[scale=0.14]
\pgfmathsetmacro\fsize{87};
\pgfmathsetmacro\xorg{20};
\pgfmathsetmacro\yorg{20};
\pgfmathsetmacro\xdlcorner{-11};
\pgfmathsetmacro\ydlcorner{-23};
\pgfmathsetmacro\stepunit{25};

\coordinate (dlcorner) at (\xdlcorner,\ydlcorner);

\draw[black,line width=0.4mm] (dlcorner) -- +(0,\fsize) -- +(\fsize,\fsize) -- +(\fsize,0) -- +(0,0);
\draw[dotted,black,line width=0.4mm] (\xdlcorner+2.5,\yorg) -- (\xdlcorner+\fsize,\yorg);
\draw[dotted,black,line width=0.4mm] (\xorg,\ydlcorner+1.8) -- (\xorg,\ydlcorner+\fsize);

\begin{scope}[every node/.style={font={\scriptsize}}]
\node at (32,\ydlcorner+1){\textbf{ Factor 1 (25.03\%)}};
\node at (\xorg-\stepunit,\ydlcorner+1){\textbf{-1}};
\node at (\xorg,\ydlcorner+1){\textbf{0}};
\node at (\xorg+\stepunit,\ydlcorner+1){\textbf{1}};
\node at (\xorg+2*\stepunit,\ydlcorner+1){\textbf{2}};
\node[rotate=90] at (\xdlcorner+1.7,7){\textbf{ Factor 2 (17.91\%)}};
\node at (-9.6,\yorg){\textbf{0}};
\node at (-9.4,\yorg-\stepunit){\textbf{-1}};
\node at (-9.6,\yorg+\stepunit){\textbf{1}};
\end{scope}

\path plot[mark=*,mark size=4mm] coordinates {(50.375,17.7) (17.55,16.5) (38.05,29.475) (36.975,14.2) (48,27.3) (3.75,40.325) (22.725,30.85) (43.425,9.675) (26.3,25.2) (11.325,34.15) (21.5,28.625) (42.35,23) (52.35,24.175) (21.1,17.45) (11.675,21.475) (13.125,30.05) (26.825,15.95) (27.525,26.425) (27.925,45.7) (61.05,45.15) (61.225,41.05) (59.175,40.275) (50.325,26.45) (45.675,23.1) (40.3,19.1) (23.9,33.1) (39.675,11.275) (5.875,39.75) (25.175,36.025) (55.55,16.05) (36.725,50.2) (44.175,15.075) (7.1,33.65) (11.575,33.15) (37.125,30.225) (37.75,40.6) (26.55,12) (41.025,18.75) (22.475,17.95) (33.275,35.625) (42.35,36.95) (48.3,13.925) (14.175,15.7) (4.475,36.775) (26.05,11.25) (49.275,12.275) (39.775,5.925) (26.175,19.975) (52.4,16.6) (50.475,26.55) (44.65,18.575) (22.65,37.05) (47.725,19.125) (31.725,24.325) (8.275,27.575) (51.9,26.275) (61.425,26.25) (-0.225,32.95) (15.65,7.825) (31.525,17.675) (39.65,10.6) (40.175,9.95) (12.55,32.85) (30.125,19.625) (33.55,9.125) (18.5,5.825) (72.575,40.525) (56.15,26.175) (16.075,10.075) (16.275,11.175) (11.45,25.375) (11.05,12.45) (17.425,-0.725) (31.725,45.35) (47.675,36.775) (43.6,24.9) (22.725,-2.05) (38.65,34.075) (8.775,38.025) (33.075,19.4) (35.6,18.525) (17.1,13) (26.175,14.275) (42.05,33.775) (7.6,28.075) (20.475,11.2) (56.075,37.15) (48.9,20.3) (-6.85,36.875) (42.925,30.3) (50.725,23.225) (34.225,43.65) (2.425,35.7) (44.425,40.975) (18.6,31.075) (26.875,22.65) (19.1,28.525) (36.9,24.075) (9.05,7.975) (37.85,15.275) (10.025,-12) (39.525,18.975) (24.725,25.025) (7.075,16.425) (9.125,26.875) (21.2,15.275) (18.775,21.075) (48.95,15.9) (45.525,38.075) (34.125,16.85) (18.625,37.225) (18.525,39.225) (37.925,43.825) (39.825,28.25) (21,24.2) (46.425,22.525) (29.2,33.425) (34.7,23.375) (56.675,33.9) (19.75,38.05) (15.025,12.525) (16.15,2.025) (1.675,31.5) (15.55,19.65) (26.7,9.8) (42.3,24.65) (6.775,28.775) (28.275,16.475) (17.075,27.075) (30.05,25.375) (27.15,19.25) (37,29.25) (2.225,38.925) (32.725,21.325) (26.1,16.425) (12.3,26) (12.25,28.525) (49.525,22.575) (41.625,31.4) (23.175,23.2) (-4.5,36.85) (8.55,-22.125) (29.75,21) (17.3,6.975) (24.075,25.05) (20.95,23.7) (20,22.65) (9.025,27.725) (-0.5,33.45) (10.575,39.95) (8.45,12.025) (7.95,17.225) (23.625,1.675) (24.65,10.575) (10.6,23.5) (5.975,16.65) (64.95,38.1) (12.5,26.675) (30.025,21.925) (31.425,15.975) (49.6,7.975) (28.925,19.1) (31.3,17.35) (12.325,36.575) (41.95,40.975) (34.275,10.525) (-2.475,35.925) (3.725,31.55) (37.525,18.9) (43.9,17.45) (11.25,36.875) (37.225,28.55) (44.075,2.925) (8.175,27.675) (49.075,29.7) (15.4,18.175) (30.45,13.325) (14,28.85) (9.725,-20.2) (21.625,27.8) (20.925,20.225) (25.575,17.625) (17.675,10.825) (18.875,30.325) (26.925,21.425) (16.45,21.375) (18.35,19.375) (29.1,31.4) (27.325,26.4) (16.3,32.35) (32.025,25.8) (34.675,31.975) (27.875,7.625) (22.45,26.55) (30.05,29.375) (29.175,30.075) (15.825,21.075) (33.125,23.275) (32.825,25.475) (25.075,6.175) (31.55,20.55) (30.75,28.375) (13.975,26.1) (24.65,21.675) (37.375,30.725) (15.15,16.9) (30.2,25.175) (24.725,17.85) (24.65,16.85) (20.2,24.2) (42.35,34.775) (16.3,15.775) (28.1,23.975) (25.925,14.425) (18.525,21.425) (31.6,12.7) (19.425,21.975)};

\begin{scope}[every node/.style={font={\tiny},inner sep=1pt,fill=white}]
\node[draw] (Kadran) at (18.95,15.075){{Kadran}};
\node[draw] (BSG3143) at (28.55,26.925){{BSG3143}};
\node[draw] (BNF1339) at (41.025,26.95){{BNF1339}};
\node[draw] (Nicolas-Chuquet) at (13.675,26.7){{Nicolas Chuquet}};
\node[draw] (BNF25) at (15.525,5.5){{BNF2050}};
\node[draw] (BNF10259) at (25.825,30.975){{BNF1259}};
\node[draw] (Traicte-praticque) at (9.175,25.45){{Traicte praticque}};
\node[draw] (Nantes456) at (29.825,14.125){{Nantes456}};
\end{scope}

\node (trouver)at (8.175,27.675){};
\node (demander)at (11.575,33.15){};
\node (item)at (17.1,13){};
\node (requerir)at (-0.5,33.45){};
\node (convenir)at (5.875,39.75){};
\node (faire)at (11.45,25.375){};
\node (vouloir)at (14,28.85){};
\node (montrer)at (13.975,26.1){};
\node (digit)at (49.275,12.275){};
\node (article)at (52.35,24.175){};
\node (nombre-compo)at (48.95,15.9){};
\node (figurer)at (25.075,6.175){};
\node (dire)at (27.875,7.675){};
\node (raison)at (16.25,15.775){};
\node (nombrateur)at (20.75,15.5){};

\begin{scope}[every node/.style={font={\footnotesize}}]
\draw[black,<-,line width=0.2mm] (trouver) -- +(-4,0)node[left]{trouver};
\draw[black,<-,line width=0.2mm] (demander) -- +(-4,-3)node[left]{demander};
\draw[black,<-,line width=0.2mm] (item) -- +(3,0)node[right]{item};
\draw[black,<-,line width=0.2mm] (requerir) -- +(-2,0)node[left]{requerir};
\draw[black,<-,line width=0.2mm] (convenir) -- +(0,3)node[above]{convenir};
\draw[black,<-,line width=0.2mm] (faire) -- +(-3,-2)node[left]{faire};
\draw[black,<-,line width=0.2mm] (vouloir) -- +(4,1)node[right]{vouloir};
\draw[black,<-,line width=0.2mm] (montrer) -- +(2.5,-2.2)node[below]{montrer};
\draw[black,<-,line width=0.2mm] (digit) -- +(4,-1)node[right]{digit};
\draw[black,<-,line width=0.2mm] (article) -- +(4,0)node[right]{article};
\draw[black,<-,line width=0.2mm] (nombre-compo) -- +(4,-1)node[right]{nombre composé};
\draw[black,<-,line width=0.2mm] (figurer) -- +(4,-3)node[right]{figurer};
\draw[black,<-,line width=0.2mm] (dire) -- +(4,-1)node[right]{dire};
\draw[black,<-,line width=0.2mm] (raison) -- +(-4,-2)node[left]{raison};
\draw[black,<-,line width=0.2mm] (nombrateur) -- +(-4,3)node[left]{nombrateur};
\end{scope}

\end{tikzpicture}
\caption{Projection on the first two factors of the FCA. The eight texts appear in frames, a few words are displayed while the remaining are simply figured by dots, for the sake of readability.}

\label{fig:CA12}
\end{figure}

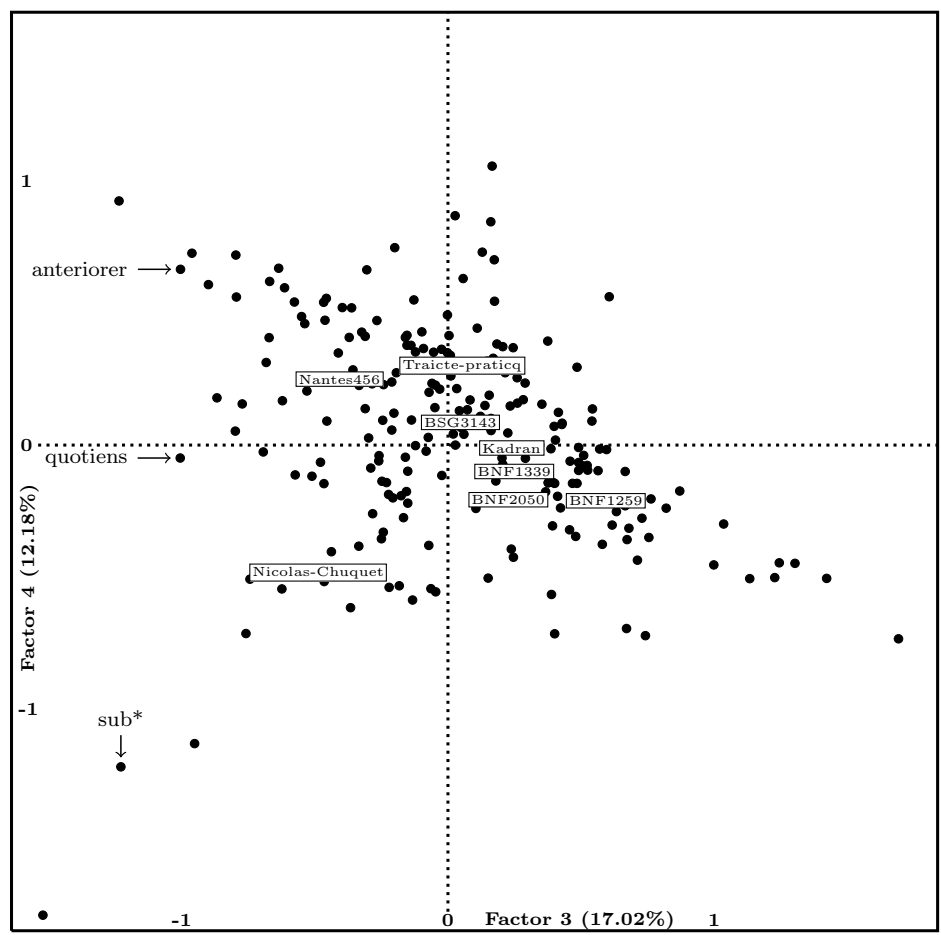
\begin{figure}[h!]
\centering

\begin{tikzpicture}[scale=0.14]

\pgfmathsetmacro\fsize{87};
\pgfmathsetmacro\xorg{20};
\pgfmathsetmacro\yorg{20};
\pgfmathsetmacro\xdlcorner{-21};
\pgfmathsetmacro\ydlcorner{-26};
\pgfmathsetmacro\stepunit{25};

\coordinate (dlcorner) at (\xdlcorner,\ydlcorner);

\draw[black,line width=0.4mm] (dlcorner) -- +(0,\fsize) -- +(\fsize,\fsize) -- +(\fsize,0) -- +(0,0);
\draw[dotted,black,line width=0.4mm] (\xdlcorner+2.5,\yorg) -- (\xdlcorner+\fsize,\yorg);
\draw[dotted,black,line width=0.4mm] (\xorg,\ydlcorner+1.8) -- (\xorg,\ydlcorner+\fsize);

\begin{scope}[every node/.style={font={\scriptsize}}]
\node at (32,\ydlcorner+1){\textbf{ Factor 3 (17.02\%)}};
\node at (\xorg-\stepunit,\ydlcorner+1){\textbf{-1}};
\node at (\xorg,\ydlcorner+1){\textbf{0}};
\node at (\xorg+\stepunit,\ydlcorner+1){\textbf{1}};
\node[rotate=90] at (\xdlcorner+1.7,7){\textbf{ Factor 4 (12.18\%)}};
\node at (-19.6,\yorg){\textbf{0}};
\node at (-19.4,\yorg-\stepunit){\textbf{-1}};
\node at (-19.6,\yorg+\stepunit){\textbf{1}};
\end{scope}

\path plot[mark=*,mark size=4mm] coordinates { (5.575,33.55) (5.65,17.175) (29.675,19.65) (16.95,28.825) (15.6,15.2) (32.275,18.35) (-18.05,-24.525) (25.075,18.775) (-5.125,36.65) (4.45,24.2) (24.35,37.55) (16.225,17.525) (32.675,18.025) (30.025,16.375) (23.775,7.4) (29.375,29.85) (16.225,14.5) (44.975,8.65) (20.7,20) (52.6,8.8) (48.35,7.35) (33.125,17.625) (36.65,17.5) (32.125,16.375) (25.15,18.175) (39.075,14.9) (9.05,9.9) (15.825,13.125) (22.75,31.075) (-2.5,35.2) (36.65,14.25) (23.875,24.725) (-3.8,-8.275) (6.55,31.5) (55.575,7.375) (0.075,38) (24.25,28.225) (24.05,21.375) (32.75,19.025) (34.25,19.625) (15.425,6.675) (16,30.2) (8.325,33.525) (30.3,15.15) (32.275,17.6) (51.125,8.85) (25.95,10.15) (8.45,31.825) (34.9,19.575) (35.15,34.05) (27.25,25.875) (9.7,28.725) (11.075,27.125) (17.95,19.425) (29.175,15.575) (6.25,32.175) (21.825,23.35) (12.55,20.675) (33.525,22.275) (19.425,17.125) (13.525,19) (12.225,23.45) (33.575,23.425) (12.225,30.3) (34.1,17.575) (16.8,33.75) (18.15,20.725) (18.575,22.225) (15,38.7) (29.725,5.85) (3.2,30.175) (4.65,34.9) (3.25,35.5) (19.95,28.75) (26.525,23.975) (27.075,24.3) (14.725,25.975) (38.225,13.075) (38.875,11.25) (19.225,25.3) (20.8,27.9) (23.05,22.725) (28.825,23.875) (22.075,24.275) (26.15,9.375) (24.6,29.575) (36.825,11.05) (38.55,1.95) (16,18.85) (62.325,1.65) (30.725,22.075) (30.1,20.475) (29.825,12.35) (12.925,13.5) (20.1,30.375) (24.025,41.15) (20.15,27.35) (13.8,16.575) (16.675,5.325) (30.575,14.05) (32.275,19.775) (8.625,22.275) (8.375,7.075) (29.425,16.45) (6.75,25.125) (20.675,41.725) (4.4,6.375) (11.9,30.7) (31.475,18.475) (32,11.35) (50.7,7.45) (19.95,32.325) (13.325,31.8) (21.5,21.025) (2.65,19.35) (14.475,6.525) (14.225,16.45) (33.075,18.075) (20.275,26.55) (4.1,36.75) (31.425,11.975) (8.575,33.9) (30.7,21.95) (25.625,21.15) (21.425,35.775) (20.825,25.35) (45.9,12.525) (8.35,16.35) (10.725,30.2) (30.375,23.1) (29.975,21.775) (16.15,29.45) (21.075,23.25) (18.65,28.8) (41.775,15.65) (15.15,26.85) (-1.7,24.475) (10.075,33.025) (20.5,21.05) (29.425,16.425) (24.475,19.975) (23.475,23.75) (18.825,25.675) (13.875,22.35) (16.575,22.375) (10.95,33) (19.4,29.075) (10.85,4.6) (29.875,16.45) (8.025,18.375) (13.925,11.75) (17.7,29.15) (14.425,15.325) (22.35,27.25) (18.4,6.4) (20.225,28.5) (0.675,23.9) (7.225,17.05) (26.5,26.375) (25.15,29.325) (17.55,30.725) (1.4,7.3) (2.925,27.825) (22.625,14) (14.925,23.025) (0.025,21.325) (21.8,22.45) (31.7,16.375) (-5.15,18.775) (40.5,14.025) (18.85,6.1) (24.15,46.425) (36.775,2.625) (12.375,36.6) (25.375,26.85) (30.025,2.125) (25.85,23.7) (11.625,10.425) (24.375,33.625) (32.125,27.375) (13.75,11.125) (12.875,25.8) (14.825,15) (27.275,18.75) (18.2,10.5) (34.5,10.6) (37,12.125) (16.525,29.45) (26.125,29.225) (35.425,12.425) (24.025,21.625) (21.225,21.675) (16.1,15.6) (18.5,25.85) (11.65,25.65) (13.95,25.725) (18.775,23.55) (-4.05,38.175) (29.375,17.275) (16.15,30.4) (13.5,18.5) (-10.725,-10.475) (23.675,27.95) (0.125,34.025) (23.225,38.275) (18.225,25) (24.5,16.6) (37.8,9.1) (1.025,2.15) (35.825,13.7) (-10.9,43.125) (12.75,17.825) (14.725,21.425) (16.95,19.95) (24.05,22.525) (28.3,19.725) (20.275,22.05) };

\begin{scope}[every node/.style={font={\tiny},inner sep=1pt,fill=white}]
\node[draw] (Kadran)at (26,19.725){{Kadran}};
\node[draw] (BSG3143)at (21.175,22.15){{BSG3143}};
\node[draw] (BNF1339)at (26.225,17.5){{BNF1339}};
\node[draw] (Nicolas-Chuquet)at (7.8,7.925){{Nicolas-Chuquet}};
\node[draw] (BNF25)at (25.675,14.825){{BNF2050}};
\node[draw] (BNF10259)at (34.825,14.725){{BNF1259}};
\node[draw] (Traicte-praticq)at (21.35,27.525){{Traicte-praticq}};
\node[draw] (Nantes456)at (9.825,26.25){{Nantes456}};
\end{scope}

\node (anteriorer)at (-5.125,36.65){};
\node (quotiens)at (-5.15,18.775){};
\node (sub)at (-10.725,-10.475){};

\begin{scope}[every node/.style={font={\footnotesize}}]
\draw[black,<-,line width=0.2mm] (sub) -- +(0,3)node[above]{sub*};
\draw[black,<-,line width=0.2mm] (quotiens) -- +(-4,0)node[left]{quotiens};
\draw[black,<-,line width=0.2mm] (anteriorer) -- +(-4,0)node[left]{anteriorer};
\end{scope}

\end{tikzpicture}

\caption{Projection on third and fourth factors of the FCA.}
\label{fig:CA34}
\end{figure}

In our case, the rows are the words ($I=219$) and the columns are the manuscripts ($J=8$). Figures~\ref{fig:CA12} and~\ref{fig:CA34} show the projection of the data on the first four factorial axes. 

The first two factors (43.94\% of the total variance) show the diversity of the cultural heritages which have built  the language of these treaties. The first factor (25.03\%) discriminates between the university legacy on the right, and the tradition of mathematical problems on the left.\\

On the left, we can observe a group whose strong homogeneity comes from its orientation towards mathematical problems (\textit{trouver} that is to say "to find", \textit{demander} which we can translate as "to ask") and their iteration (\textit{item, idem}). That vocabulary can be found most often in both the appendix of \textit{Triparty en la science des nombres} (Nicolas Chuquet) and \textit{Le Traicte de la praticque}. Furthermore, there are more verbal forms on this side of the axis than on the other. And we can find verbs like \textit{requerir} which means "to require", \textit{convenir} "to agree", \textit{faire} "to do". Some of them are prescriptive, as \textit{devoir} "to have to" or \textit{vouloir} "to want" for example, while others introduce examples, as \textit{montrer} "to show". All these texts contain a lot of mathematical problems and in a way are practical. On the right, the texts of BnF. fr. 1339 and Med. Nantes 456 are clearly more representative of the university culture, containing latin words sequences.\\

The second axis (17.91\% of the variance) is mostly characterized by the manuscript of BNF. fr. 2050 and also by \textit{Kadran aux marchans}. It displays words of Italo-Provencal origin, like \textit{nombrateur} which refers to the division's numerator. Designations of the fraction and operation of division take a significant part of the information while the most contributory words (for ex. \textit{figurer} "to draw") allow us to examine another dimension of these works: the graphical representation as a continuation of writing.\\

The following factors 3, 4 and further show the Lexicon that seems to be more related to the singularity of some manuscripts. The importance of Nicolas Chuquet inertia of  factors 3 and 4 singles out this book on the plane (see Figure~\ref{fig:CA34}) in relation to the rest of the corpus.\\

With the manuscript of Nantes 456, at left, factor 3 highlights a vocabulary of some technical accuracy, in any case rare in the rest of the corpus, like \textit{quotiens} "quotient", \textit{anteriorer} "to put before". At right, there is a very diversified vocabulary, associated to manuscrit 10259, which is a well organized compilation of a copy of \textit{Kadrans aus marchans} and of a lot of problems whose origin has not been fully identified.\\ 

Correspondence Analysis displays the particularities of each text, but leaves untouched some more complex elements of the data. For instance, we have to see the third factor to understand that the  \textit{Triparty en la science des nombres} (Nicolas Chuquet) and the \textit{Traicte de la praticque} use different university mathematical cultures. These two treaties are not only copying university algorithms as they were taught at university at that time, they have their own originality.\\

Moreover, we cannot assert that the words which appear in the center of the graph represent a 'common vocabulary': as a matter of fact, we should  analyze all the successive factors in order to build the list of words constituting the 'common vocabulary'. It is a very cumbersome task.\\

\subsection{Kohonen Maps}

SOM-based algorithms were very often used for text mining purpose. Oja and Kaski's seminal book \cite{oja99} provides a lot of examples on this field. A major tool for that purpose is WEBSOM method and software\footnote{See \url{http://websom.hut.fi/websom/}}, as defined for instance in \cite{kohonen99,lagus99,lagus04}. Other important papers (among hundreds) are \cite{rauber03,yang04,honkela}.

Most of them look for classification and clustering using keywords, put in evidence the main features, associate documents with their most characteristic words, to define proximity in order to define clusters and hierarchies between documents. Techniques such as WEBSOM are especially designed to deal with massive documents collections.

Our purpose is very different, since we have very few documents and since we look for a subset of words which are not "specific" of some manuscript, but contrarily belong to a common vocabulary.\\

Factorial Correspondence Analysis (FCA) suffers from some limitations as explained in section\ref{fcadef}. To overcome this, we use a variant of the SOM algorithm which deals with the same kind of data, i. e. a contingency table. This variant of SOM was defined in  \cite{cottrell98} or \cite{cottrell05} and we refer to it as  KORRESP algorithm. Let us recall this definition.\\

Tha data are displayed as explained in fourth paragraph of section \ref{fcadef} in a contingency $I=219$ by $J=8$ table. The data are normalized applying equation (\ref{normalization}), exactly in the same way as for Factorial Correspondence Analysis. The normalized contingency table is denoted by 
$T^{norm}$ where:
$$t^{norm}_{i,j}= \frac{t_{i,j}}{\sqrt{\sum_i t_{i,j} \sum_j t_{i,j}}}.$$

We consider a Kohonen map, and associate to each unit $u$ a code-vector $C_u$ with $(J+I)$ components. The first $J$ components evolve in the space of the rows (the words), while the last $I$ components belong to the space of the columns (the manuscripts).

Let us denote 
\begin{equation}
C_u = (C_J, C_I)_u = (C_{J,u}, C_{I,u}),
\end{equation}
to put in evidence the structure of the code-vector $C_u$.\\

We use the SOM algorithm as a double learning process, by alternatively drawing a $T^{norm}$ row (a word) and a $T^{norm}$ column (a manuscript).

When we draw a row $i$, we associate the column $j(i)$ that maximizes the coefficient $t^{norm}_{i,j}$, so: 
\begin{equation}
j(i) = \arg\max_{j} t^{norm}_{i,j} = \arg\max_{j}\frac{t_{i,j}}{\sqrt{\sum_i t_{i,j} \sum_j t_{i,j}}}
\end{equation}

that maximizes the conditional probability of $j$ given $i$.
We then create an extended $(J+I)$~-~dimensional row vector $X = (i, j(i)) = (X_J, X_I)$. See Figure~\ref{fig:latefig}.

\begin{center}
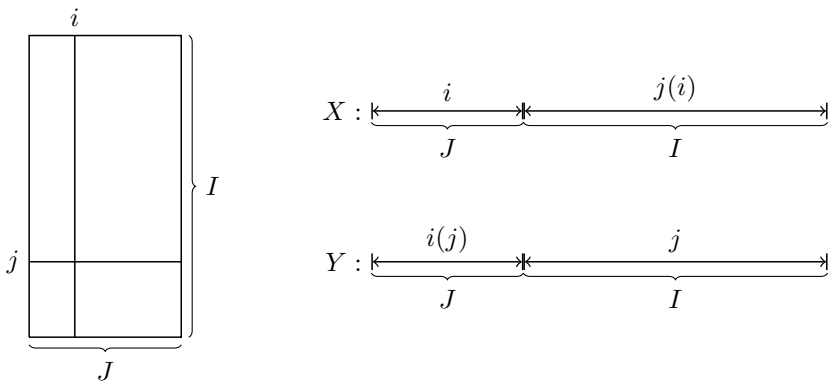
\begin{figure}
	\begin{tikzpicture}
		\draw[black,line width=0.2mm] (1,0) -- +(0,4) -- +(2,4) -- +(2,0) -- +(0,0);
		\draw[black,line width=0.2mm] (1,1)node[left]{$j$} -- +(2,0);
		\draw[black,line width=0.2mm] (1.6,0) -- +(0,4)node[above]{$i$};
		\draw[black,decorate,decoration={brace,raise = 0.1cm}] (3,4) -- +(0,-4)node[right=0.2cm,pos=0.5]{$I$};
		\draw[black,decorate,decoration={brace,raise = 0.1cm}] (3,0) -- +(-2,0)node[below=0.2cm,pos=0.5]{$J$};
		\draw[black,line width=0.2mm,|<->|] (7.5,1) -- +(4,0)node[above,pos=0.5]{$j$};
		\draw[black,line width=0.2mm,|<->|] (5.5,1)node[left]{$Y:$} -- +(2,0)node[above,pos=0.5]{$i(j)$};
		\draw[black,line width=0.2mm,|<->|] (7.5,3) -- +(4,0)node[above,pos=0.5]{$j(i)$};
		\draw[black,line width=0.2mm,|<->|] (5.5,3)node[left]{$X:$} -- +(2,0)node[above,pos=0.5]{$i$};
		\draw[black,decorate,decoration={brace,raise = 0.15cm}] (7.5,1) -- +(-2,0)node[below=0.25cm,pos=0.5]{$J$};
		\draw[black,decorate,decoration={brace,raise = 0.15cm}] (11.5,1) -- +(-4,0)node[below=0.25cm,pos=0.5]{$I$};
		\draw[black,decorate,decoration={brace,raise = 0.15cm}] (7.5,3) -- +(-2,0)node[below=0.25cm,pos=0.5]{$J$};
		\draw[black,decorate,decoration={brace,raise = 0.15cm}] (11.5,3) -- +(-4,0)node[below=0.25cm,pos=0.5]{$I$};
	\end{tikzpicture}
\label{fig:latefig}
\caption{Building of the extended, symmetrized table in the KORRESP algorithm}
\end{figure}
\end{center}

Subsequently, we look for the closest of all the code vectors, in terms of the Euclidean distance restricted to the first $J$ components. Note $u_0$ the winning unit. Next we move the code-vector of the unit $u_0$ and its neighbors towards the extended vector $X = (i, j(i))$, as per the customary Kohonen law. Let us write down the formal definition:

\begin{eqnarray}
u_0 &=& \arg\min_u \|X_J - C_{J,u}\|\\
C^{new}_u &=& C^{old}_u + \epsilon \sigma (u,u_0) (X -C^{old}_u) 
\end{eqnarray}

where $\epsilon$ is the adaptation parameter (positive, decreasing with time), and $\sigma$ is the neighborhood function, such that 
$\sigma(u,u_0) = 1$ if $u$ and $u_0$ are neighbour in the Kohonen network, and $\sigma(u,u_0) = 0$ if not.

The reason to associate a row and a column in such a way is to keep the row-column associations which are realized in classical FCA by the fact that the principal axes of both Principal Component Analysis are strongly correlated.

The procedure is the same when we draw a column $j$ with dimension $I$ (a column of $T^{norm}$). We associate the row $i(j)$ that maximizes the coefficient $t^{norm}_{i,j}$, so: 
\begin{equation}
i(j) = \arg\max_{i} t^{norm}_{i,j} = \arg\max_{i}\frac{t_{i,j}}{\sqrt{\sum_i t_{i,j} \sum_j t_{i,j}}}
\end{equation}

that maximizes the conditional probability of $i$ given $j$.
We then create an extended $(J+I)$-dimensional column vector 
$Y = (i(j), j) = (Y_J,Y_I)$. 

We then seek the code-vector that is the closest, in terms of the Euclidean distance restricted to the last $I$ components. Let $v_0$ be the winning unit. Next we move the code-vector of the unit $v_0$ and its neighbors towards the extended vector $Y = (i(j), j)$, as per the customary Kohonen law. Let us write down the formal definition:

\begin{eqnarray}
v_0 &=& \arg\min_v \|Y_I - C_{I,u}\| \\
C^{new}_v &=& C^{old}_v + \epsilon \sigma (v,v_0) (Y - C^{old}_v) 
\end{eqnarray}
 
where $\epsilon$ and $\sigma$ are defined as before.

This two-steps computation carries out a Kohonen classification of the rows (the words), together with a classification of the columns, maintaining all the while the associations of both rows and columns.

We can sum up the definition of the KORRESP algorithm 
\begin{itemize}
\item \textit{normalization} of the rows and of the columns in the way as in FCA computation,
\item \textit{definition} of an extended data table by associating to each row the most probable column and to each column the most probable row, 
\item  \textit{simultaneous classification} of the rows and of the columns onto a Kohonen map, by using the rows of the extended data table as input for the SOM algorithm.
\end{itemize}

After convergence of the training step, the items of the rows and of the columns are simultaneously classified. In our example, one can see proximity between words, between texts, between words and texts. It is the same goal as in Factorial Correspondence Analysis. The advantage is that it is not necessary to examine several projection planes: the whole information can be read on the Kohonen Map.\\

We display below (Figure~\ref{fig:koho1}) the SOM map which simultaneously represents the words and the texts. For this map as for all the remaining of this paper we use the online algorithm, a 10$\times$10 grid and the following simple neighborhood function: 1 for the eight (fewer if we are along one edge of the map) nodes adjacent to the selected one and 0 for the others.\\

One can observe that the interpretation (see Figure \ref{fig:robsteph}) is very similar to the interpretation that could be done from the Factorial Correspondence Analysis projections.
But, as an example of the robustness problem, we can compare two different Kohonen maps (in Figures~\ref{fig:koho1} and~\ref{fig:koho2}) and the respective positions of the words \textit{raison} "reason" and \textit{dire} "to say", very far from each other in the first map while neighboring in the second one.\\

\hyphenpenalty=10000
\setlength{\tabcolsep}{1mm}
\setlength{\arrayrulewidth}{1pt}

\begin{figure}[ht!]
\centering
\scalebox{0.6}{
\begin{tabular}{|p{1.9cm}|p{1.9cm}|p{1.9cm}|p{1.9cm}|p{1.9cm}|p{1.9cm}|p{1.9cm}|p{1.9cm}|p{1.9cm}|p{1.9cm}|}
\hline
minutes super* & & notes & calculer cubic & fois mettre \textbf{BNF10259} & & {contraire} {depenser} {falloir} {meme} {racine} & aller donc ensuivre {savoir} & {multiplier} & {regle} venir\\
\hline
gecter & & & & & dessous & barrater {demi} & somme & voir & assembler\\
\hline
 & \textbf{BSG3143} & parteur & defaillir duplation mediation nommer numeration senestre & semblables & & circulaires demeurer derenier disaine ecrire entendre formes nombrateur oter prouver & entrer laisser rien & emprunter figure regarder & figure de~non~rien fraction muer rayes retenir\\
\hline
notables & & denomi- nations multipli- cateur nominateur ordonne & abaisser comptes endroit proposer & anteriorer diminution enseignement enseigner moyen signifiant surplus trancher & possible {reduire} & difficile progression repondre & avaluer & faillir & monter partiteur\\
\hline
bref gectons {multipli-}{cation} {pratique} {seulement} & generale latin {nombrer} proportion reduction unite & soustraction & denomi- nateur & entier & ajoutement & & remotion \textbf{Kadran} & & \textbf{BNF2050}\\
\hline
chose& nulle & {ensemble} {partie} & \textbf{Nantes456} & & & & {partement} {valoir} & & etre\\
\hline
arithmetique {compter} preuve tenir & cubbement destre digit diviser diviseur division lignes nombre compo{se}querir & \cellcolor{grisclair} \textbf{\textcolor{red}{DIRE}} {figurer} poser {science} & & & abreger lever precedent quotiens & numerateur & moindre nombre \textbf{Nicolas Chuquet} & partir plus \cellcolor{grisclair} \textbf{\textcolor{red}{RAISON}} & \\
\hline
egalir egaliser especes question total traiter & mesurer & grand & apparaitre {position} {soustraire} & {ajouter} & leurs prendre quant quantefois reponse trouver & & {part} & commun Item & devoir {droit} {exemple} reste rester\\
\hline
algorisme article carrees cercle envient ligne pair sain \textbf{BNF1339} & addition chiffre former & {double} {doubler} moitie & appeler {donner} {maniere} {pouvoir} {se}& bailler demander mises nomper pareillement vouloir & & egale faire montrer necessaire romp selon & & & \\
\hline
cautelle & demontrer dessus & {garder} {regle} {de} {trois} & aliquot composer corps moins sub* toutefois & partant plaisir proportionel- lement requerir residu & appartenir convenir demande difference egaulx maieur millions rate survendre tant & & naturel roupt \textbf{Traicte praticque} & & fausse\\
\hline
\end{tabular}
}
\caption{Example of Kohonen Map. Manuscripts are in bold. Notice that \textit{raison} (9,7) and \textit{dire} (3,7) are far apart from each other.}
\label{fig:koho1}
\end{figure}

\begin{figure}[ht!]
\begin{center}
\begin{tikzpicture}
\pgfmathsetmacro\fsize{10};

\draw[black,line width=0.4mm] (0,0) -- +(0,\fsize) -- +(\fsize,\fsize) -- +(\fsize,0) -- +(0,0);
\draw[dotted,black,line width=0.4mm] (0,0) -- (\fsize,\fsize);

\begin{scope}[shape=rectangle,font={\small}]
\node[draw] at (0.9,2.2){BNF 1339};
\node[draw] at (1.3,9){BSG 3143};
\node[draw] at (4,6){Nantes 456};
\node[draw] at (5,9.5){BNF 10259};
\node[draw] at (8,1){Traicté de la praticque};
\node[draw] at (8,5.2){Nicolas Chuquet};
\node[draw] at (8,6.2){Kadran aus Marchans};
\node[draw] at (9.4,6.8){BNF 2050};
\end{scope}

\begin{scope}[every node/.style={font={\scriptsize}}]
\node[text width=2cm] at (0.6,0)[below left]{closer to university algorithms};
\node[text width=3cm] at (1,\fsize)[above left]{influenced by \textit{l'algorisme à jetons} and the manuscript from Jean Adam};
\node[text width=3.5cm] at (\fsize-1,-0.2)[below right]{influenced by the \textit{Traité de la praticque}, attempting to make mathematical knowledge more accessible to a merchant audience};
\node[text width=2.5cm] at (\fsize-1,\fsize)[above right]{closer to a practical vocabulary};
\end{scope}

\end{tikzpicture}
\caption{Interpretation of the Kohonen map (Figure~\ref{fig:koho1}), the diagonal opposes university and practical vocabularies}
\label{fig:robsteph}

\end{center}
\end{figure}
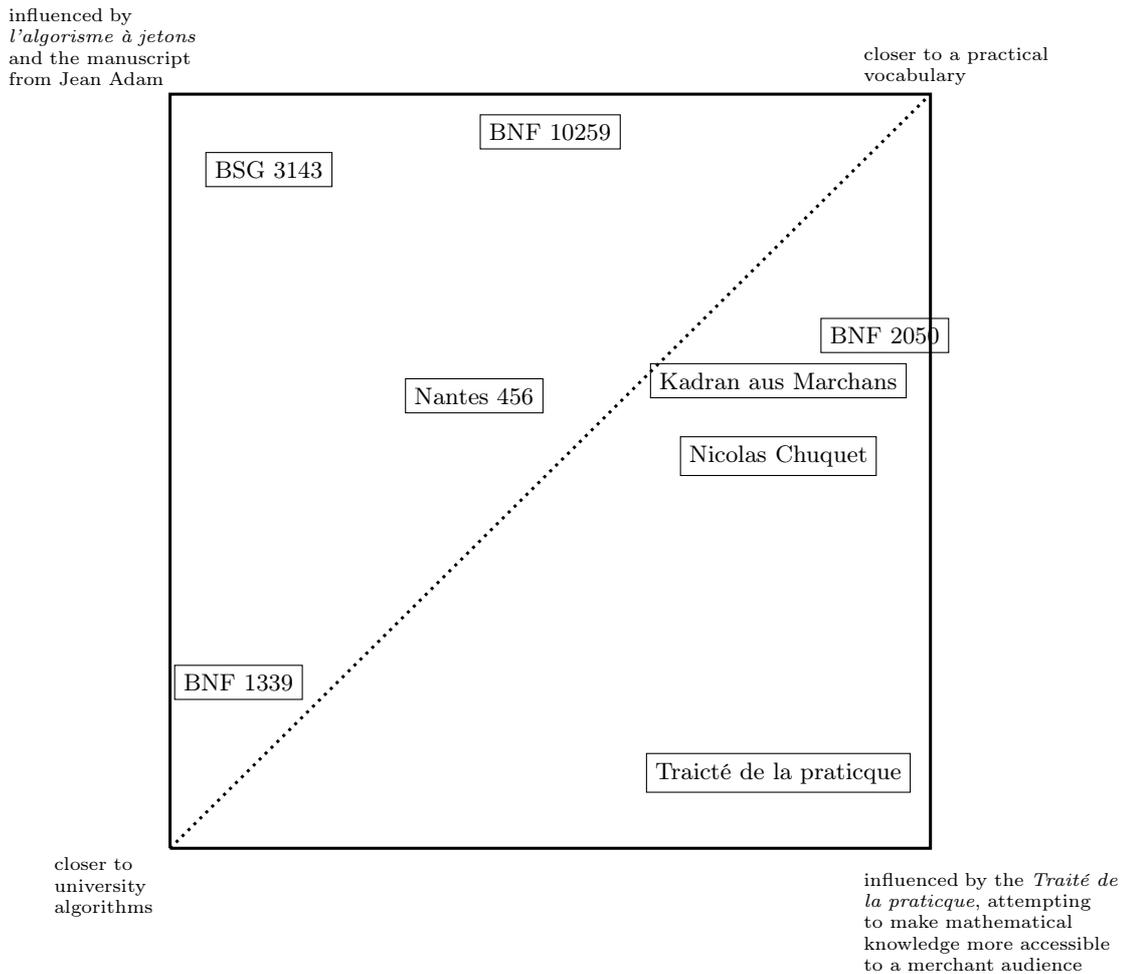

\hyphenpenalty=10000
\setlength{\tabcolsep}{1mm}
\setlength{\arrayrulewidth}{1pt}

\begin{figure}[ht!]
\centering
\scalebox{0.6}{
\begin{tabular}{|p{1.95cm}|p{1.9cm}|p{1.95cm}|p{1.88cm}|p{1.88cm}|p{1.88cm}|p{1.88cm}|p{1.88cm}|p{1.88cm}|p{1.88cm}|}
\hline
 cautelle latin nombrer pair & demontrer & dessus \textbf{BNF1339} & addition traiter & compter double doubler somme & & devoir droit multiplier rien & barrater contraire depenser falloir meme racine & & calculer cubic fois mettre \textbf{BNF10259} \\ \hline
 unite & article chiffre former & & &demeurer preuve & remotion & & & dessous & \\ \hline
 duplation mediation & &carrees cercle envient ligne sain & &progression tenir & & faillir & &formes oter & demi \\ \hline
 anteriorer cubbement destre digit diviser diviseur division lignes moyen nombrecompose signifiant & egalir especes querir & algorisme egaliser total & arithmetique circulaires mesurer question & garder laisser & & & & disaine ecrire entendre nombrateur possible & prouver reduire \\ \hline 
 ajoutement apparaitre diminution ensemble partie precedent proposer semblables & science soustraction & \cellcolor{grisclair} \textbf{\textcolor{red}{DIRE}} figurer poser & ajouter moitie position pouvoir \cellcolor{grisclair} \textbf{\textcolor{red}{RAISON}} se soustraire & moins regle de trois & partement & fausse & exemple & avaluer derenier difficile entrer repondre Kadran & \\ \hline
\textbf{Nantes456} & & & assembler grand valoir & appeler composer donner & aller bailler mises & convenir demande demander donc nomper pareillement & partant & & \\ \hline
 & voir & pratique regle & chose & & & ensuivre plaisir requerir residu tant & faire & monter & emprunter figure figure de non rien muer partiteur rayes regarder retenir \textbf{BNF2050} \\ \hline
 denominateur entier & fraction & bref generale multiplication & nulle proportion & part seulement & montrer vouloir & leurs necessaire quant savoir selon trouver & rester venir & reste & \\ \hline
 comptes endroit enseignement enseigner nommer senestre surplus trancher & parteur & nominateur notes & gectons ordonne & egale & romp roupt & & & abreger commun item lever nombre partir plus prendre reponse & \\ \hline
 abaisser defaillir denominations multiplicateur numeration reduction & & gecter minutes super  \textbf{BSG3143} & notables & appartenir difference egaulx maieur millions rate survendre & naturel & etre  \textbf{Traicte praticque} & & maniere moindre numerateur quantefois quotiens & aliquot corps proportionel- lement sub toutefois \textbf{Nicolas Chuquet}\\
\hline
\end{tabular}
}
\caption{Another example of Kohonen Map. This time, \textit{raison} (4,5) and \textit{dire} (3,5) are neighbors.}
\label{fig:koho2}
\end{figure}

The Kohonen algorithm is stochastic, and it can  happen that several runs get different results, and that these differences can be troublesome. Hence the idea to introduce repetitions of the runs to separate stable and robust results from purely stochastic behavior. In the following, we study the variability of the maps which provides new information.

\section{Getting extra information through the extraction of fickle words}\label{method}

In its classical presentation \cite{kohonen95,cottrell98}, the SOM algorithm is an iterative algorithm, which takes as input a dataset $\mathbf{x}_i, i\in\{1,\hdots,N\}$ and computes code-vectors $\mathbf{m}_u, u\in\{1,\hdots,U\}$ which define the map.\\

We know that self-organization is reached at the end of the algorithm, which implies that close data in the input space have to belong to the same class or to neighboring classes, that is to say that they are projected on the same prototypes or on neighboring prototypes on the map. In what follows, we call neighbors data that belong either to the same unit or to two adjacent units. But the reciprocal is not exact: for a given run of the algorithm, two given data can be neighbors on the map, while they are not in the input space. That drawback comes from the fact that there is no perfect fit between a two-dimensional map and the data space (except when the intrinsic dimension is exactly 2). As we just notice, since the SOM algorithm is a stochastic one, the resulting maps can be different from one run to another. How to overcome this difficulty?\\

In fact, we can use this drawback to improve the interpretation and the analysis of relations between the studied words. Our hypothesis is that the repetitive use of this method can help us to identify words that are strongly attracted/repulsed and also fickle pairs.

\subsection{Neighborhood and robustness of information on Kohonen maps}

We address the issue of computing a reliability level for the neighboring (or no-neighboring) relations in a SOM map. More precisely, if we consider several runs of the SOM algorithm, for a given size of the map and for a given data set, we observe that most of pairs are almost always neighbors or always not neighbors. But there are also pairs whose associations look random. These pairs are called \textit{fickle} pairs. This question was addressed by \cite{debodt02} in a bootstrap frame.\\

According to their paper, we can define:
$NEIGH_{i,j}^l=0$ if $x_i$ and $x_j$ are not neighbors in the $l$-th run of the algorithm, and 
$NEIGH_{i,j}^l=1$ if $x_i$ and $x_j$ are neighbors in the $l$-th run of the algorithm, 
where $(x_i, x_j)$ is a given pair of data, $l$ is the number of the observed runs of the SOM algorithm.\\

Then $Y_{i,j}= \sum_{l=1}^L NEIGH_{i,j}^l$ is the number of times when the data $x_i$ and $x_j$ are neighbor for $L$ different, independent runs. 
The stability index $\mathcal{M}_{i,j}$ is defined as the average of $NEIGH_{i,j}$ over all the runs ($l=1, \hdots, L$), i. e. 

\begin{equation}
\mathcal{M}_{i,j} = \frac{\sum_{l=1}^L NEIGH_{i,j}^l}{L}=\frac{Y_{i,j}}{L}.
\end{equation}
 
The next step is to compare it to the value it would have if the data $x_i$ and $x_j$ were neighbors by chance in a completely random way.\\

So we can use a classical statistical test to check the significance of the stability index $\mathcal{M}_{i,j}$. Let $U$ be the number of units on the map. If edge effects are not taken into account, the number of units involved in a neighborhood region (as defined here) is 9 in a two-dimensional map. So for a fixed pair of data $x_i$ and $x_j$, the probability of being neighbors in a random way is equal to $9/U$ (it is the probability for $x_j$ to be a neighbor of $x_i$ by chance once the class $x_i$ belongs to is determined).\\

As $Y_{i,j}= \sum_{l=1}^L NEIGH_{i,j}^l$ is the number of times when the data $x_i$ and $x_j$ are neighbor for $L$ different, independent runs, it is easy to see that $Y_{i,j}$ is distributed as a Binomial distribution with parameters $L$ and $9/U$. 

Using the classical approximation of Binomial Distribution by a Gaussian one ($L$ is large and $9/U$ not too small), we can build the critical region of the test of null hypothesis $H_0$ "$x_i$ and $x_j$ are neighbors by chance" against hypothesis $H_1$: " the fact that $x_i$ and $x_j$ are neighbors or not is significant".\\

We conclude that the critical region for a test level of $5\%$ based on $Y_{i,j}$, is 
\begin{equation} 
]-\infty,  L \frac{9}{U} - 1.96 \sqrt{L \frac{9}{U} ( 1 - \frac{9}{U})}[ 
\quad   \bigcup   \quad  ] L \frac{9}{U} + 1.96 \sqrt{L \frac{9}{U} ( 1 - \frac{9}{U}) }, +\infty [ 
\end{equation}
 
For the frequency (i.e. the stability index) $\mathcal{M}_{i,j} = Y_{i,j}/L$,  the critical region is 

\begin{equation} 
]-\infty,  \frac{9}{U} - 1.96 \sqrt{ \frac{9}{U L} ( 1 - \frac{9}{U})}[ 
\quad   \bigcup   \quad  ]  \frac{9}{U} + 1.96 \sqrt{  \frac{9}{U L} ( 1 - \frac{9}{U}) }, +\infty [ 
\end{equation}

To simplify the notations, , let us put
\begin{equation}\label{AB}
 A = \frac{9}{U} \text{ and } B = 1.96 \sqrt{\frac{9}{UL} ( 1 - \frac{9}{U})}.
\end{equation}

Then, practically, for each pair of words $(x_i,x_j)$, we  compute the index $\mathcal{M}_{i,j} = Y_{i,j}/L$, and apply the following rule:\\ 
 
\begin{itemize}
\item if their index is greater than $A + B$, they are almost always neighbors in a significant way, the words attract each other.
\item if their index is comprised between $A - B$ and $A+B$, their proximity is due to randomness, they are a fickle pair.
\item if their index is less than $A - B$, they are almost never neighbor, the words repulse each other.
\end{itemize}

\subsection{Identification of fickle pairs}

We run KORRESP $L$ times and store the result in a matrix $\mathcal{M}$ of size $(N+p)\times(N+p)$. The value stored in a given cell $i,j$ is the proportion of maps where $i$ and $j$ are neighbors.\\

Table~\ref{bertin1} displays an example of the  first nine rows and columns of such a matrix. We have highlighted with colors three different situations. According to the theoretical study mentioned above:

\begin{itemize}
	\item Black cells stand for pairs that are neighbors with high probability (proximity happens with frequency greater than $A+B$, here $0.1787$).
	\item White cells stand for pairs that are not neighbors with high probability (proximity happens with frequency less than $A-B$, here $0.0014$).
	\item Grey cells are not conclusive, they are the \textit{fickle pairs}.
\end{itemize}

If we rearrange the order of cells and columns through Bertin permutations, we immediately make remarkable clustering properties appear (see Table~\ref{bertin2}.)\\

For each word, through this treatment we get a list of words that can roughly be grouped around two poles: the strongly associated and the almost never associated ones. Between these two extremes lies a central yet difficult to characterize.

This technique could be used for classification, but here our main objective is a bit different: we are mostly interested in a characterization of words that have high mobility in Kohonen maps, that we call \textit{fickle words}.

\label{m1}

\begin{table}
\begin{small}
\begin{tabular}{|l|l|l|l|l|l|l|l|l|l|}
\hline
 & abaisser & abreger & addition & ajoutement & ajouter & algorisme & aliquot & aller & anteriorer \\ \hline
abaisser  & \cellcolor{gris} 1 & 0 & \cellcolor{grisclair} 0.025 & \cellcolor{gris} 0.275 & 0 & \cellcolor{grisclair} 0.05 & 0 & 0 & \cellcolor{gris} 0.525\\ \hline
abreger & 0 & \cellcolor{gris} 1 & 0 & 0 & \cellcolor{gris} 0.25 & 0 & \cellcolor{gris} 0.325 & 0 & \cellcolor{grisclair} 0.025\\ \hline
addition & \cellcolor{grisclair} 0.025 & 0 & \cellcolor{gris} 1 & 0 & 0 & \cellcolor{gris} 0.875 & 0 & \cellcolor{grisclair} 0.05 & 0\\ \hline
ajoutement  & \cellcolor{gris} 0.275 & 0 & 0 & \cellcolor{gris} 1 & \cellcolor{grisclair} 0.025 & 0 & 0 & \cellcolor{grisclair} 0.025 & \cellcolor{gris} 0.7\\ \hline
ajouter & 0 & \cellcolor{gris} 0.25 & 0 & \cellcolor{grisclair} 0.025 & \cellcolor{gris} 1 & \cellcolor{grisclair} 0.025 & \cellcolor{grisclair} 0.15 & \cellcolor{grisclair} 0.125 & 0\\ \hline
algorisme & \cellcolor{grisclair} 0.05 & 0 & \cellcolor{gris} 0.875 & 0 & \cellcolor{grisclair} 0.025 & \cellcolor{gris} 1 & 0 & 0 & 0\\ \hline
aliquot & 0 & \cellcolor{gris} 0.325 & 0 & 0 & \cellcolor{grisclair} 0.15 & 0 & \cellcolor{gris} 1 & \cellcolor{grisclair} 0.025 & 0\\ \hline
aller & 0 & 0 & \cellcolor{grisclair} 0.05 & \cellcolor{grisclair} 0.025 & \cellcolor{grisclair} 0.125 & 0 & \cellcolor{grisclair} 0.025 & \cellcolor{gris} 1 & 0\\ \hline
anteriorer & \cellcolor{gris} 0.525 & \cellcolor{grisclair} 0.025 & 0 & \cellcolor{gris} 0.7 & 0 & 0 & 0 & 0 & \cellcolor{gris} 1 \\ \hline
\end{tabular}
\end{small}
\caption{Frequency of neighborhood matrix (excerpt)}
\label{bertin1}
\end{table}

\begin{table}
\begin{small}

\begin{tabular}{|l|l|l|l|l|l|l|l|l|l|}
\hline
 & abaisser & ajoutement & anteriorer& abreger & ajouter & aliquot & addition & algorisme & aller \\ \hline
abaisser  & \cellcolor{gris} 1 & \cellcolor{gris} 0.275 & \cellcolor{gris} 0.525 & 0 & 0 & 0 & \cellcolor{grisclair} 0.025 & \cellcolor{grisclair} 0.05 & 0 \\ \hline
ajoutement  & \cellcolor{gris} 0.275 & \cellcolor{gris} 1 & \cellcolor{gris} 0.7 & 0 & \cellcolor{grisclair} 0.025 & 0 & 0 & 0 & \cellcolor{grisclair} 0.025 \\ \hline
anteriorer & \cellcolor{gris} 0.525 & \cellcolor{gris} 0.7 & \cellcolor{gris} 1 & \cellcolor{grisclair} 0.025 & 0 & 0 & 0 & 0 & 0 \\ \hline
abreger & 0 & 0 &\cellcolor{grisclair} 0.025& \cellcolor{gris} 1 & \cellcolor{gris} 0.25 & \cellcolor{gris} 0.325 & 0 & 0 & 0 \\ \hline
ajouter & 0 & \cellcolor{grisclair} 0.025 & 0& \cellcolor{gris} 0.25 & \cellcolor{gris} 1 & \cellcolor{grisclair} 0.15 & 0 & \cellcolor{grisclair} 0.025 & \cellcolor{grisclair} 0.125 \\ \hline
aliquot & 0 & 0 & 0& \cellcolor{gris} 0.325 & \cellcolor{grisclair} 0.15 & \cellcolor{gris} 1 & 0 & 0 & \cellcolor{grisclair} 0.025 \\ \hline
addition & \cellcolor{grisclair} 0.025 & 0 & 0& 0 & 0 & 0 & \cellcolor{gris} 1 & \cellcolor{gris} 0.875 & \cellcolor{grisclair} 0.05 \\ \hline
algorisme & \cellcolor{grisclair} 0.05 & 0 & 0& 0 & \cellcolor{grisclair} 0.025 & 0 & \cellcolor{gris} 0.875 & \cellcolor{gris} 1 & 0 \\ \hline
aller & 0 & \cellcolor{grisclair} 0.025 & 0& 0 & \cellcolor{grisclair} 0.125 & \cellcolor{grisclair} 0.025 & \cellcolor{grisclair} 0.05 & 0 & \cellcolor{gris} 1 \\ \hline
\end{tabular}
\end{small}
\caption{Frequency of neighborhood matrix (same excerpt as \ref{bertin1}, with row and columns reorganized)}
\label{bertin2}
\end{table}

\subsection{From fickle pairs to fickle words}
\label{fromfrom}

We call \textit{fickle} a word which belongs to a huge number of fickle pairs:

$$|\{i,|\mathcal{M}_{i,j}-A| \leq B \}| \geq \Theta$$

Unfortunately, it is not quite an easy task to find an appropriate threshold $\Theta$. Here we have decided to fix it according to data interpretation. The $30$ ficklest words, whose number of safe neighbors/non-neighbors (non-fickle pairs) is between $89$ and $119$, are displayed in Figure~\ref{fig:ficklest}.

\begin{center}
\begin{figure}
\begin{small}
\begin{tabular}{|lll|}
\hline
\textit{contraire} "opposite" (89) & \textit{regle de trois} "rule of three" (104) & \textit{depenser} "to expend" (112) \\
\textit{doubler} "to double" (89) & \textit{savoir}  "to know" (105) & \textit{racine} "root" (113) \\
\textit{falloir} "to need" (93) & \textit{partie} "to divide" (105) & \textit{chose} "thing" (113) \\
\textit{meme} "same, identical" (93) &  \textit{position} "position" (107) & \textit{compter} "to count" (113) \\
\textit{pratique} "practical" (94) &  \textit{exemple} "for example" (107) & \textit{dire} "to say" (113) \\
\textit{seulement}  "only" (94) &  \textit{demi} "half"  (108) & \textit{nombrer} "count" (115) \\
\textit{double} "double" (97) &  \textit{garder} "to keep"(109) & \textit{raison} "calculation, problem" (116) \\
\textit{multiplication} (99) &  \textit{science} "science" (109) & \textit{donner} "to give" (117) \\
\textit{reduire} "to reduce" (103) & \textit{pouvoir} "can" (111) & \textit{ensemble} "together" (117) \\
\textit{regle} "rule" (103) & \textit{se} "if" (111) & \textit{valoir} "to be worth" (119)\\
\hline
\end{tabular}

\end{small}
\caption{30 ficklest words among 219 studied. For each word, the number between brackets stands for how many non-fickle pairs it belongs to.}
\label{fig:ficklest}
\end{figure}
\end{center}

\subsection{Graph of robust neighborhood}

Let us have a different look at the neighborhood matrix $(f_{ij})$ where $f_{ij}$ is the frequency of two words  belonging to the same neighborhood. Instead of trying to jump right ahead and identify fickle words in an absolute way, we can study the robust connections between words \textit{per se}, in order to produce some interesting clustering of the words.\\

For example, if we have a look at the excerpt from Table~\ref{bertin1}, we notice immediately that some groups of words are very often in the same neighborhood, while their connections to the rest of the graph are much more hazardous. This initial intuition becomes quite obvious if we reorganize the rows and columns (following Bertin's permutation matrices idea), as we can see on Table~\ref{bertin2}.

We cannot display here the whole matrix for the 219 forms - in addition, the algorithm for reorganization would not be efficient enough - so we have decided to focus on a specific group of words: the fickle words. Indeed, the fickle words are the most difficult to study, since by definition they do not have a very fixed position on the Kohonen maps, and additionally it appears that they are not well distinguished by Factorial Correspondence Analysis either.

\begin{table}[h!]
\centering
\setlength{\tabcolsep}{1pt}
\setlength{\arrayrulewidth}{0.5pt}
\scalebox{0.72}{ 
\begin{tabular}{|p{2pt}l|*{30}{c|}}
\hline
 && \begin{sideways}contraire\end{sideways} & \begin{sideways}doubler\end{sideways} & \begin{sideways}falloir\end{sideways} & \begin{sideways}meme\end{sideways} & \begin{sideways}pratique\end{sideways} & \begin{sideways}seulement\end{sideways} & \begin{sideways}double\end{sideways} & \begin{sideways}multiplication \end{sideways} & \begin{sideways}reduire\end{sideways} & \begin{sideways}regle\end{sideways} & \begin{sideways}regle de trois\end{sideways} & \begin{sideways}partie\end{sideways} & \begin{sideways}savoir\end{sideways} & \begin{sideways}exemple\end{sideways} & \begin{sideways}position\end{sideways} & \begin{sideways}demi\end{sideways} & \begin{sideways}garder\end{sideways} & \begin{sideways}science\end{sideways} & \begin{sideways}pouvoir\end{sideways} & \begin{sideways}se\end{sideways} & \begin{sideways}depenser\end{sideways} & \begin{sideways}chose\end{sideways} & \begin{sideways}compter\end{sideways} & \begin{sideways}dire\end{sideways} & \begin{sideways}racine\end{sideways} & \begin{sideways}nombrer\end{sideways} & \begin{sideways}raison\end{sideways} & \begin{sideways}donner\end{sideways} & \begin{sideways}ensemble\end{sideways} & \begin{sideways}valoir\end{sideways}
\\ \hline
&contraire & \cellcolor{gris}1 & 0 & \cellcolor{gris}.6{\tiny 25} & \cellcolor{gris}.7{\tiny 5} & 0 & .0{\tiny 25} & 0 & .0{\tiny 75} & \cellcolor{gris}.2 & .1{\tiny 25} & .0{\tiny 5} & .0{\tiny 5} & \cellcolor{gris}.4{\tiny 25} & .0{\tiny 5} & .0{\tiny 25} & \cellcolor{gris}.6{\tiny 5} & 0 & 0 & .0{\tiny 5} & .0{\tiny 25} & \cellcolor{gris}.5{\tiny 5} & .1 & 0 & 0 & \cellcolor{gris}.6{\tiny 75} & .0{\tiny 25} & .0{\tiny 5} & .0{\tiny 75} & .0{\tiny 75} & .0{\tiny 75}
\\ \hline
&doubler & 0 & \cellcolor{gris}1 & .0{\tiny 5} & .0{\tiny 25} & .0{\tiny 5} & .0{\tiny 75} & \cellcolor{gris}.9{\tiny 25} & .0{\tiny 75} & .0{\tiny 5} & .0{\tiny 5} & .0{\tiny 75} & \cellcolor{gris}.5{\tiny 75} & .0{\tiny 5} & .0{\tiny 25} & .1{\tiny 75} & .0{\tiny 25} & \cellcolor{gris}.2{\tiny 25} & \cellcolor{gris}.3{\tiny 25} & .0{\tiny 75} & \cellcolor{gris}.2{\tiny 5} & .0{\tiny 75} & \cellcolor{gris}.2{\tiny 75} & \cellcolor{gris}.6{\tiny 5} & .1 & .0{\tiny 25} & .1 & .0{\tiny 25} & \cellcolor{gris}.4{\tiny 5} & \cellcolor{gris}.3{\tiny 75} & \cellcolor{gris}.2
\\ \hline
&falloir & \cellcolor{gris}.6{\tiny 25} & .0{\tiny 5} & \cellcolor{gris}1 & \cellcolor{gris}.7{\tiny 75} & .0{\tiny 25} & 0 & .0{\tiny 5} & 0 & \cellcolor{gris}.2 & .1 & .0{\tiny 5} & .1{\tiny 25} & \cellcolor{gris}.2{\tiny 5} & .0{\tiny 5} & .1 & \cellcolor{gris}.5{\tiny 5} & .1 & .0{\tiny 5} & .1 & .0{\tiny 25} & \cellcolor{gris}.8{\tiny 5} & .1 & .0{\tiny 75} & .1{\tiny 25} & \cellcolor{gris}.6{\tiny 75} & .0{\tiny 25} & .1{\tiny 5} & .1{\tiny 75} & .1{\tiny 5} & .0{\tiny 5}
\\ \hline
&meme & \cellcolor{gris}.7{\tiny 5} & .0{\tiny 25} & \cellcolor{gris}.7{\tiny 75} & \cellcolor{gris}1 & 0 & .0{\tiny 25} & .0{\tiny 25} & .0{\tiny 5} & \cellcolor{gris}.2{\tiny 75} & .1 & .0{\tiny 5} & .1{\tiny 25} & \cellcolor{gris}.2{\tiny 75} & .0{\tiny 5} & .0{\tiny 75} & \cellcolor{gris}.7{\tiny 25} & .0{\tiny 5} & .0{\tiny 25} & .0{\tiny 75} & .0{\tiny 25} & \cellcolor{gris}.7{\tiny 5} & .1{\tiny 75} & .0{\tiny 5} & .0{\tiny 75} & \cellcolor{gris}.6{\tiny 75} & .0{\tiny 25} & .1{\tiny 25} & .1{\tiny 5} & .1{\tiny 5} & .0{\tiny 5}
\\ \hline
&pratique & 0 & .0{\tiny 5} & .0{\tiny 25} & 0 & \cellcolor{gris}1 & \cellcolor{gris}.7{\tiny 75} & .0{\tiny 5} & \cellcolor{gris}.6{\tiny 25} & \cellcolor{gris}.2{\tiny 75} & \cellcolor{gris}.3{\tiny 25} & .1 & .1{\tiny 5} & .1 & \cellcolor{gris}.3{\tiny 75} & .0{\tiny 5} & 0 & .0{\tiny 75} & .0{\tiny 5} & .0{\tiny 5} & .1{\tiny 5} & .0{\tiny 5} & \cellcolor{gris}.4{\tiny 25} & .1 & 0 & 0 & .0{\tiny 5} & .0{\tiny 25} & .1{\tiny 25} & .1 & .1
\\ \hline
&seulement & .0{\tiny 25} & .0{\tiny 75} & 0 & .0{\tiny 25} & \cellcolor{gris}.7{\tiny 75} & \cellcolor{gris}1 & .0{\tiny 75} & \cellcolor{gris}.6{\tiny 75} & \cellcolor{gris}.2{\tiny 5} & \cellcolor{gris}.2{\tiny 75} & .1{\tiny 25} & .1{\tiny 5} & .1{\tiny 25} & \cellcolor{gris}.3{\tiny 5} & .0{\tiny 5} & 0 & .0{\tiny 5} & .0{\tiny 75} & .0{\tiny 5} & .1{\tiny 75} & .0{\tiny 25} & \cellcolor{gris}.4{\tiny 75} & .1 & 0 & 0 & .0{\tiny 5} & .0{\tiny 25} & .1{\tiny 75} & .0{\tiny 75} & .0{\tiny 25}
\\ \hline
&double & 0 & \cellcolor{gris}.9{\tiny 25} & .0{\tiny 5} & .0{\tiny 25} & .0{\tiny 5} & .0{\tiny 75} & \cellcolor{gris}1 & .0{\tiny 25} & .0{\tiny 75} & .0{\tiny 5} & .1{\tiny 75} & \cellcolor{gris}.5 & .0{\tiny 5} & .0{\tiny 5} & .1{\tiny 75} & .0{\tiny 25} & \cellcolor{gris}.3{\tiny 5} & \cellcolor{gris}.2{\tiny 75} & .1 & \cellcolor{gris}.2{\tiny 75} & .1 & \cellcolor{gris}.2{\tiny 75} & \cellcolor{gris}.6{\tiny 25} & .1 & .0{\tiny 25} & .0{\tiny 75} & .0{\tiny 5} & \cellcolor{gris}.4{\tiny 75} & \cellcolor{gris}.3{\tiny 25} & \cellcolor{gris}.2
\\ \hline
&multiplication & .0{\tiny 75} & .0{\tiny 75} & 0 & .0{\tiny 5} & \cellcolor{gris}.6{\tiny 25} & \cellcolor{gris}.6{\tiny 75} & .0{\tiny 25} & \cellcolor{gris}1 & .1 & \cellcolor{gris}.3{\tiny 25} & .1 & .1{\tiny 25} & .0{\tiny 5} & .0{\tiny 5} & 0 & .0{\tiny 5} & .0{\tiny 25} & .0{\tiny 75} & .0{\tiny 25} & .0{\tiny 75} & 0 & \cellcolor{gris}.5 & .0{\tiny 25} & .0{\tiny 25} & 0 & \cellcolor{gris}.3{\tiny 25} & .0{\tiny 25} & .0{\tiny 25} & .0{\tiny 5} & 0
\\ \hline
&reduire & \cellcolor{gris}.2 & .0{\tiny 5} & \cellcolor{gris}.2 & \cellcolor{gris}.2{\tiny 75} & \cellcolor{gris}.2{\tiny 75} & \cellcolor{gris}.2{\tiny 5} & .0{\tiny 75} & .1 & \cellcolor{gris}1 & \cellcolor{gris}.2{\tiny 75} & \cellcolor{gris}.2{\tiny 25} & \cellcolor{gris}.2{\tiny 25} & .0{\tiny 75} & \cellcolor{gris}.2{\tiny 75} & .1 & \cellcolor{gris}.3{\tiny 25} & \cellcolor{gris}.3{\tiny 25} & .1 & 0 & .1{\tiny 25} & .1{\tiny 5} & \cellcolor{gris}.3{\tiny 75} & .1 & .0{\tiny 5} & .1{\tiny 5} & .0{\tiny 75} & .0{\tiny 5} & .0{\tiny 75} & \cellcolor{gris}.2{\tiny 75} & .0{\tiny 5}
\\ \hline
&regle & .1{\tiny 25} & .0{\tiny 5} & .1 & .1 & \cellcolor{gris}.3{\tiny 25} & \cellcolor{gris}.2{\tiny 75} & .0{\tiny 5} & \cellcolor{gris}.3{\tiny 25} & \cellcolor{gris}.2{\tiny 75} & \cellcolor{gris}1 & .0{\tiny 25} & \cellcolor{gris}.2{\tiny 25} & .1{\tiny 75} & .1 & 0 & .0{\tiny 5} & .0{\tiny 5} & \cellcolor{gris}.2{\tiny 5} & .0{\tiny 5} & .0{\tiny 5} & .0{\tiny 75} & \cellcolor{gris}.2{\tiny 25} & .1 & \cellcolor{gris}.3 & .0{\tiny 5} & .1{\tiny 5} & 0 & .0{\tiny 25} & \cellcolor{gris}.2{\tiny 5} & .1
\\ \hline
&regle de trois & .0{\tiny 5} & .0{\tiny 75} & .0{\tiny 5} & .0{\tiny 5} & .1 & .1{\tiny 25} & .1{\tiny 75} & .1 & \cellcolor{gris}.2{\tiny 25} & .0{\tiny 25} & \cellcolor{gris}1 & .0{\tiny 25} & 0 & .1{\tiny 5} & \cellcolor{gris}.7 & .0{\tiny 5} & \cellcolor{gris}.6{\tiny 5} & .0{\tiny 25} & \cellcolor{gris}.5{\tiny 5} & \cellcolor{gris}.6{\tiny 25} & .0{\tiny 25} & \cellcolor{gris}.2{\tiny 25} & .0{\tiny 25} & .0{\tiny 5} & .0{\tiny 25} & .0{\tiny 75} & \cellcolor{gris}.6{\tiny 75} & \cellcolor{gris}.4{\tiny 75} & 0 & \cellcolor{gris}.3
\\ \hline
&partie & .0{\tiny 5} & \cellcolor{gris}.5{\tiny 75} & .1{\tiny 25} & .1{\tiny 25} & .1{\tiny 5} & .1{\tiny 5} & \cellcolor{gris}.5 & .1{\tiny 25} & \cellcolor{gris}.2{\tiny 25} & \cellcolor{gris}.2{\tiny 25} & .0{\tiny 25} & \cellcolor{gris}1 & .0{\tiny 75} & .0{\tiny 5} & .1 & .0{\tiny 5} & \cellcolor{gris}.2 & \cellcolor{gris}.5{\tiny 5} & .0{\tiny 5} & .1{\tiny 75} & .1{\tiny 5} & \cellcolor{gris}.3{\tiny 75} & \cellcolor{gris}.6{\tiny 5} & \cellcolor{gris}.2{\tiny 75} & .0{\tiny 25} & .1 & 0 & \cellcolor{gris}.3 & \cellcolor{gris}.7{\tiny 75} & .1{\tiny 25}
\\ \hline
&savoir & \cellcolor{gris}.4{\tiny 25} & .0{\tiny 5} & \cellcolor{gris}.2{\tiny 5} & \cellcolor{gris}.2{\tiny 75} & .1 & .1{\tiny 25} & .0{\tiny 5} & .0{\tiny 5} & .0{\tiny 75} & .1{\tiny 75} & 0 & .0{\tiny 75} & \cellcolor{gris}1 & .1{\tiny 25} & .0{\tiny 75} & \cellcolor{gris}.3 & 0 & .0{\tiny 25} & .1 & .0{\tiny 5} & .1{\tiny 75} & .0{\tiny 75} & .1 & .0{\tiny 25} & \cellcolor{gris}.2{\tiny 5} & 0 & .1 & \cellcolor{gris}.2 & 0 & \cellcolor{gris}.2
\\ \hline
&exemple & .0{\tiny 5} & .0{\tiny 25} & .0{\tiny 5} & .0{\tiny 5} & \cellcolor{gris}.3{\tiny 75} & \cellcolor{gris}.3{\tiny 5} & .0{\tiny 5} & .0{\tiny 5} & \cellcolor{gris}.2{\tiny 75} & .1 & .1{\tiny 5} & .0{\tiny 5} & .1{\tiny 25} & \cellcolor{gris}1 & .0{\tiny 75} & .0{\tiny 75} & .1 & .0{\tiny 5} & .0{\tiny 75} & \cellcolor{gris}.3 & .0{\tiny 5} & .1{\tiny 5} & .1 & .0{\tiny 25} & .0{\tiny 25} & 0 & .1{\tiny 5} & \cellcolor{gris}.3 & .0{\tiny 25} & \cellcolor{gris}.3{\tiny 25}
\\ \hline
&position & .0{\tiny 25} & .1{\tiny 75} & .1 & .0{\tiny 75} & .0{\tiny 5} & .0{\tiny 5} & .1{\tiny 75} & 0 & .1 & 0 & \cellcolor{gris}.7 & .1 & .0{\tiny 75} & .0{\tiny 75} & \cellcolor{gris}1 & .0{\tiny 75} & \cellcolor{gris}.6 & .1{\tiny 25} & \cellcolor{gris}.7{\tiny 75} & \cellcolor{gris}.7{\tiny 25} & .0{\tiny 5} & .1 & .0{\tiny 5} & \cellcolor{gris}.2 & 0 & 0 & \cellcolor{gris}.7{\tiny 25} & \cellcolor{gris}.5{\tiny 25} & .0{\tiny 75} & \cellcolor{gris}.5{\tiny 25}
\\ \hline
&demi & \cellcolor{gris}.6{\tiny 5} & .0{\tiny 25} & \cellcolor{gris}.5{\tiny 5} & \cellcolor{gris}.7{\tiny 25} & 0 & 0 & .0{\tiny 25} & .0{\tiny 5} & \cellcolor{gris}.3{\tiny 25} & .0{\tiny 5} & .0{\tiny 5} & .0{\tiny 5} & \cellcolor{gris}.3 & .0{\tiny 75} & .0{\tiny 75} & \cellcolor{gris}1 & .0{\tiny 75} & 0 & .0{\tiny 75} & .0{\tiny 5} & \cellcolor{gris}.5{\tiny 75} & .0{\tiny 75} & 0 & .0{\tiny 25} & \cellcolor{gris}.6{\tiny 25} & .0{\tiny 25} & .1{\tiny 75} & .1{\tiny 5} & .0{\tiny 75} & .0{\tiny 5}
\\ \hline
&garder & 0 & \cellcolor{gris}.2{\tiny 25} & .1 & .0{\tiny 5} & .0{\tiny 75} & .0{\tiny 5} & \cellcolor{gris}.3{\tiny 5} & .0{\tiny 25} & \cellcolor{gris}.3{\tiny 25} & .0{\tiny 5} & \cellcolor{gris}.6{\tiny 5} & \cellcolor{gris}.2 & 0 & .1 & \cellcolor{gris}.6 & .0{\tiny 75} & \cellcolor{gris}1 & \cellcolor{gris}.2 & \cellcolor{gris}.4{\tiny 25} & \cellcolor{gris}.5 & .0{\tiny 75} & \cellcolor{gris}.2 & .1{\tiny 75} & .1 & .0{\tiny 25} & .0{\tiny 75} & \cellcolor{gris}.4{\tiny 75} & \cellcolor{gris}.3{\tiny 25} & .1{\tiny 75} & .1{\tiny 25}
\\ \hline
&science & 0 & \cellcolor{gris}.3{\tiny 25} & .0{\tiny 5} & .0{\tiny 25} & .0{\tiny 5} & .0{\tiny 75} & \cellcolor{gris}.2{\tiny 75} & .0{\tiny 75} & .1 & \cellcolor{gris}.2{\tiny 5} & .0{\tiny 25} & \cellcolor{gris}.5{\tiny 5} & .0{\tiny 25} & .0{\tiny 5} & .1{\tiny 25} & 0 & \cellcolor{gris}.2 & \cellcolor{gris}1 & .1{\tiny 25} & .1 & .0{\tiny 25} & .1{\tiny 5} & \cellcolor{gris}.4{\tiny 75} & \cellcolor{gris}.3{\tiny 5} & 0 & .1 & .0{\tiny 25} & .0{\tiny 75} & \cellcolor{gris}.6{\tiny 25} & .0{\tiny 75}
\\ \hline
&pouvoir & .0{\tiny 5} & .0{\tiny 75} & .1 & .0{\tiny 75} & .0{\tiny 5} & .0{\tiny 5} & .1 & .0{\tiny 25} & 0 & .0{\tiny 5} & \cellcolor{gris}.5{\tiny 5} & .0{\tiny 5} & .1 & .0{\tiny 75} & \cellcolor{gris}.7{\tiny 75} & .0{\tiny 75} & \cellcolor{gris}.4{\tiny 25} & .1{\tiny 25} & \cellcolor{gris}1 & \cellcolor{gris}.6{\tiny 5} & .0{\tiny 75} & .0{\tiny 25} & .0{\tiny 25} & \cellcolor{gris}.2{\tiny 25} & .0{\tiny 25} & .0{\tiny 25} & \cellcolor{gris}.8 & \cellcolor{gris}.3{\tiny 25} & .0{\tiny 5} & \cellcolor{gris}.4{\tiny 5}
\\ \hline
&se & .0{\tiny 25} & \cellcolor{gris}.2{\tiny 5} & .0{\tiny 25} & .0{\tiny 25} & .1{\tiny 5} & .1{\tiny 75} & \cellcolor{gris}.2{\tiny 75} & .0{\tiny 75} & .1{\tiny 25} & .0{\tiny 5} & \cellcolor{gris}.6{\tiny 25} & .1{\tiny 75} & .0{\tiny 5} & \cellcolor{gris}.3 & \cellcolor{gris}.7{\tiny 25} & .0{\tiny 5} & \cellcolor{gris}.5 & .1 & \cellcolor{gris}.6{\tiny 5} & \cellcolor{gris}1 & 0 & \cellcolor{gris}.2 & .1{\tiny 5} & .0{\tiny 75} & 0 & .0{\tiny 5} & \cellcolor{gris}.6{\tiny 5} & \cellcolor{gris}.6{\tiny 75} & .1{\tiny 25} & \cellcolor{gris}.5{\tiny 5}
\\ \hline
&depenser & \cellcolor{gris}.5{\tiny 5} & .0{\tiny 75} & \cellcolor{gris}.8{\tiny 5} & \cellcolor{gris}.7{\tiny 5} & .0{\tiny 5} & .0{\tiny 25} & .1 & 0 & .1{\tiny 5} & .0{\tiny 75} & .0{\tiny 25} & .1{\tiny 5} & .1{\tiny 75} & .0{\tiny 5} & .0{\tiny 5} & \cellcolor{gris}.5{\tiny 75} & .0{\tiny 75} & .0{\tiny 25} & .0{\tiny 75} & 0 & \cellcolor{gris}1 & .1{\tiny 75} & .1{\tiny 75} & .0{\tiny 25} & \cellcolor{gris}.8 & .0{\tiny 25} & .0{\tiny 75} & .1{\tiny 75} & .1 & 0
\\ \hline
&chose & .1 & \cellcolor{gris}.2{\tiny 75} & .1 & .1{\tiny 75} & \cellcolor{gris}.4{\tiny 25} & \cellcolor{gris}.4{\tiny 75} & \cellcolor{gris}.2{\tiny 75} & \cellcolor{gris}.5 & \cellcolor{gris}.3{\tiny 75} & \cellcolor{gris}.2{\tiny 25} & \cellcolor{gris}.2{\tiny 25} & \cellcolor{gris}.3{\tiny 75} & .0{\tiny 75} & .1{\tiny 5} & .1 & .0{\tiny 75} & \cellcolor{gris}.2 & .1{\tiny 5} & .0{\tiny 25} & \cellcolor{gris}.2 & .1{\tiny 75} & \cellcolor{gris}1 & \cellcolor{gris}.4{\tiny 25} & .0{\tiny 25} & .0{\tiny 75} & \cellcolor{gris}.2{\tiny 5} & .0{\tiny 25} & .1{\tiny 5} & .1{\tiny 5} & .0{\tiny 5}
\\ \hline
&compter & 0 & \cellcolor{gris}.6{\tiny 5} & .0{\tiny 75} & .0{\tiny 5} & .1 & .1 & \cellcolor{gris}.6{\tiny 25} & .0{\tiny 25} & .1 & .1 & .0{\tiny 25} & \cellcolor{gris}.6{\tiny 5} & .1 & .1 & .0{\tiny 5} & 0 & .1{\tiny 75} & \cellcolor{gris}.4{\tiny 75} & .0{\tiny 25} & .1{\tiny 5} & .1{\tiny 75} & \cellcolor{gris}.4{\tiny 25} & \cellcolor{gris}1 & .0{\tiny 75} & 0 & .0{\tiny 5} & 0 & \cellcolor{gris}.3{\tiny 75} & \cellcolor{gris}.4{\tiny 25} & .1{\tiny 25}
\\ \hline
&dire & 0 & .1 & .1{\tiny 25} & .0{\tiny 75} & 0 & 0 & .1 & .0{\tiny 25} & .0{\tiny 5} & \cellcolor{gris}.3 & .0{\tiny 5} & \cellcolor{gris}.2{\tiny 75} & .0{\tiny 25} & .0{\tiny 25} & \cellcolor{gris}.2 & .0{\tiny 25} & .1 & \cellcolor{gris}.3{\tiny 5} & \cellcolor{gris}.2{\tiny 25} & .0{\tiny 75} & .0{\tiny 25} & .0{\tiny 25} & .0{\tiny 75} & \cellcolor{gris}1 & 0 & .0{\tiny 25} & .0{\tiny 5} & .0{\tiny 25} & \cellcolor{gris}.4{\tiny 25} & .1{\tiny 25}
\\ \hline
&racine & \cellcolor{gris}.6{\tiny 75} & .0{\tiny 25} & \cellcolor{gris}.6{\tiny 75} & \cellcolor{gris}.6{\tiny 75} & 0 & 0 & .0{\tiny 25} & 0 & .1{\tiny 5} & .0{\tiny 5} & .0{\tiny 25} & .0{\tiny 25} & \cellcolor{gris}.2{\tiny 5} & .0{\tiny 25} & 0 & \cellcolor{gris}.6{\tiny 25} & .0{\tiny 25} & 0 & .0{\tiny 25} & 0 & \cellcolor{gris}.8 & .0{\tiny 75} & 0 & 0 & \cellcolor{gris}1 & .0{\tiny 75} & 0 & .0{\tiny 25} & .0{\tiny 5} & .0{\tiny 25}
\\ \hline
&nombrer & .0{\tiny 25} & .1 & .0{\tiny 25} & .0{\tiny 25} & .0{\tiny 5} & .0{\tiny 5} & .0{\tiny 75} & \cellcolor{gris}.3{\tiny 25} & .0{\tiny 75} & .1{\tiny 5} & .0{\tiny 75} & .1 & 0 & 0 & 0 & .0{\tiny 25} & .0{\tiny 75} & .1 & .0{\tiny 25} & .0{\tiny 5} & .0{\tiny 25} & \cellcolor{gris}.2{\tiny 5} & .0{\tiny 5} & .0{\tiny 25} & .0{\tiny 75} & \cellcolor{gris}1 & 0 & 0 & .1 & 0
\\ \hline
&raison & .0{\tiny 5} & .0{\tiny 25} & .1{\tiny 5} & .1{\tiny 25} & .0{\tiny 25} & .0{\tiny 25} & .0{\tiny 5} & .0{\tiny 25} & .0{\tiny 5} & 0 & \cellcolor{gris}.6{\tiny 75} & 0 & .1 & .1{\tiny 5} & \cellcolor{gris}.7{\tiny 25} & .1{\tiny 75} & \cellcolor{gris}.4{\tiny 75} & .0{\tiny 25} & \cellcolor{gris}.8 & \cellcolor{gris}.6{\tiny 5} & .0{\tiny 75} & .0{\tiny 25} & 0 & .0{\tiny 5} & 0 & 0 & \cellcolor{gris}1 & \cellcolor{gris}.3{\tiny 75} & 0 & \cellcolor{gris}.4{\tiny 25}
\\ \hline
&donner & .0{\tiny 75} & \cellcolor{gris}.4{\tiny 5} & .1{\tiny 75} & .1{\tiny 5} & .1{\tiny 25} & .1{\tiny 75} & \cellcolor{gris}.4{\tiny 75} & .0{\tiny 25} & .0{\tiny 75} & .0{\tiny 25} & \cellcolor{gris}.4{\tiny 75} & \cellcolor{gris}.3 & \cellcolor{gris}.2 & \cellcolor{gris}.3 & \cellcolor{gris}.5{\tiny 25} & .1{\tiny 5} & \cellcolor{gris}.3{\tiny 25} & .0{\tiny 75} & \cellcolor{gris}.3{\tiny 25} & \cellcolor{gris}.6{\tiny 75} & .1{\tiny 75} & .1{\tiny 5} & \cellcolor{gris}.3{\tiny 75} & .0{\tiny 25} & .0{\tiny 25} & 0 & \cellcolor{gris}.3{\tiny 75} & \cellcolor{gris}1 & .1{\tiny 75} & \cellcolor{gris}.5{\tiny 5}
\\ \hline
&ensemble & .0{\tiny 75} & \cellcolor{gris}.3{\tiny 75} & .1{\tiny 5} & .1{\tiny 5} & .1 & .0{\tiny 75} & \cellcolor{gris}.3{\tiny 25} & .0{\tiny 5} & \cellcolor{gris}.2{\tiny 75} & \cellcolor{gris}.2{\tiny 5} & 0 & \cellcolor{gris}.7{\tiny 75} & 0 & .0{\tiny 25} & .0{\tiny 75} & .0{\tiny 75} & .1{\tiny 75} & \cellcolor{gris}.6{\tiny 25} & .0{\tiny 5} & .1{\tiny 25} & .1 & .1{\tiny 5} & \cellcolor{gris}.4{\tiny 25} & \cellcolor{gris}.4{\tiny 25} & .0{\tiny 5} & .1 & 0 & .1{\tiny 75} & \cellcolor{gris}1 & .1{\tiny 25}
\\ \hline
&valoir & .0{\tiny 75} & \cellcolor{gris}.2 & .0{\tiny 5} & .0{\tiny 5} & .1 & .0{\tiny 25} & \cellcolor{gris}.2 & 0 & .0{\tiny 5} & .1 & \cellcolor{gris}.3 & .1{\tiny 25} & \cellcolor{gris}.2 & \cellcolor{gris}.3{\tiny 25} & \cellcolor{gris}.5{\tiny 25} & .0{\tiny 5} & .1{\tiny 25} & .0{\tiny 75} & \cellcolor{gris}.4{\tiny 5} & \cellcolor{gris}.5{\tiny 5} & 0 & .0{\tiny 5} & .1{\tiny 25} & .1{\tiny 25} & .0{\tiny 25} & 0 & \cellcolor{gris}.4{\tiny 25} & \cellcolor{gris}.5{\tiny 5} & .1{\tiny 25} & \cellcolor{gris}1
\\ \hline
\end{tabular}
}
\caption{Frequency of neighborhood matrix for the ficklest words only = adjacency matrix of the neighborhood graph of the ficklest}
\label{bertinfickle1}
\end{table}

Table~\ref{bertinfickle1} shows the frequency matrix for the 30 ficklest words. The clustering is not obvious \textit{a priori}, so we can use a different representation for better visualization of the underlying structures. We can fix the  threshold $A+B$ as defined in equation~(\ref{AB}) and consider this matrix as the adjacency matrix of a graph $G(V,E)$ such that:\\

\begin{itemize}
 \item the set of vertices $V$ is identified to the fickle words
 \item the set of edges $E$ is defined by $(i,j) \in E \Leftrightarrow f_{ij} > A+B$
\end{itemize}

In other terms, $G$ is the graph of highly probable neighborhood relations in Kohonen maps.

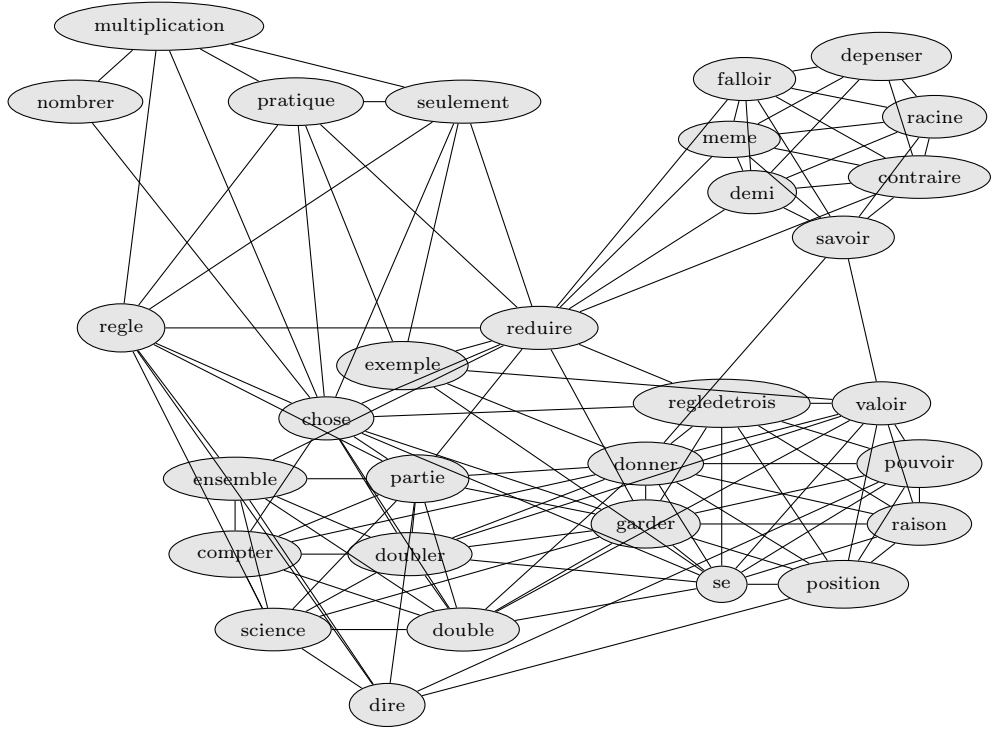
\begin{figure}[h!]
\centering
\begin{tikzpicture}[scale=1]
\begin{scope}[every node/.style={font={\scriptsize},ellipse,draw,fill=grisclair}]
\node(chose) at (3.2,4.8){chose};
\node(compter) at (2,3){compter};
\node(contraire) at (11,8){contraire};
\node(demi) at (8.8,7.8){demi};
\node(depenser) at (10.5,9.6){depenser};
\node(dire) at (4,1){dire};
\node(donner) at (7.4,4.2){donner};
\node(double) at (5,2){double};
\node(doubler) at (4.3,3){doubler};
\node(ensemble) at (2,4){ensemble};
\node(exemple) at (4.2,5.5){exemple};
\node(falloir) at (8.7,9.3){falloir};
%\node(figurer) at (1,0){figurer};
\node(garder) at (7.4,3.4){garder};
\node(meme) at (8.5,8.5){meme};
\node(multiplication) at (1,10){multiplication};
\node(nombrer) at (-0.1,9){nombrer};
\node(partie) at (4.4,4){partie};
\node(position) at (10,2.6){position};
\node(pouvoir) at (11,4.2){pouvoir};
\node(pratique) at (2.8,9){pratique};
\node(racine) at (11.2,8.8){racine};
\node(raison) at (11,3.4){raison};
\node(reduire) at (6,6){reduire};
\node(regle) at (0.5,6){regle};
\node(regledetrois) at (8.4,5){regledetrois};
\node(savoir) at (10,7.2){savoir};
\node(science) at (2.5,2){science};
\node(se) at (8.4,2.6){se};
\node(seulement) at (5,9){seulement};
\node(valoir) at (10.5,5){valoir};
\end{scope}

\draw(compter) -- (chose);
\draw(demi) -- (contraire);
\draw(depenser) -- (contraire);
\draw(depenser) -- (demi);
\draw(donner) -- (compter);
\draw(double) -- (chose);
\draw(double) -- (compter);
\draw(double) -- (donner);
\draw(doubler) -- (chose);
\draw(doubler) -- (compter);
\draw(doubler) -- (donner);
\draw(doubler) -- (double);
\draw(ensemble) -- (compter);
\draw(ensemble) -- (dire);
\draw(ensemble) -- (double);
\draw(ensemble) -- (doubler);
\draw(exemple) -- (donner);
\draw(falloir) -- (contraire);
\draw(falloir) -- (demi);
\draw(falloir) -- (depenser);
\draw(garder) -- (chose);
\draw(garder) -- (donner);
\draw(garder) -- (double);
\draw(garder) -- (doubler);
\draw(meme) -- (contraire);
\draw(meme) -- (demi);
\draw(meme) -- (depenser);
\draw(meme) -- (falloir);
\draw(multiplication) -- (chose);
\draw(nombrer) -- (chose);
\draw(nombrer) -- (multiplication);
\draw(partie) -- (chose);
\draw(partie) -- (compter);
\draw(partie) -- (dire);
\draw(partie) -- (donner);
\draw(partie) -- (double);
\draw(partie) -- (doubler);
\draw(partie) -- (ensemble);
%\draw(partie) -- (figurer);
\draw(partie) -- (garder);
\draw(position) -- (dire);
\draw(position) -- (donner);
%\draw(position) -- (figurer);
\draw(position) -- (garder);
\draw(pouvoir) -- (dire);
\draw(pouvoir) -- (donner);
%\draw(pouvoir) -- (figurer);
\draw(pouvoir) -- (garder);
\draw(pouvoir) -- (position);
\draw(pratique) -- (chose);
\draw(pratique) -- (exemple);
\draw(pratique) -- (multiplication);
\draw(racine) -- (contraire);
\draw(racine) -- (demi);
\draw(racine) -- (depenser);
\draw(racine) -- (falloir);
\draw(racine) -- (meme);
\draw(raison) -- (donner);
\draw(raison) -- (garder);
\draw(raison) -- (position);
\draw(raison) -- (pouvoir);
\draw(reduire) -- (chose);
\draw(reduire) -- (contraire);
\draw(reduire) -- (demi);
\draw(reduire) -- (ensemble);
\draw(reduire) -- (exemple);
\draw(reduire) -- (falloir);
\draw(reduire) -- (garder);
\draw(reduire) -- (meme);
\draw(reduire) -- (partie);
\draw(reduire) -- (pratique);
\draw(regle) -- (chose);
\draw(regle) -- (dire);
\draw(regle) -- (ensemble);
%\draw(regle) -- (figurer);
\draw(regle) -- (multiplication);
\draw(regle) -- (partie);
\draw(regle) -- (pratique);
\draw(regle) -- (reduire);
\draw(regledetrois) -- (chose);
\draw(regledetrois) -- (donner);
\draw(regledetrois) -- (garder);
\draw(regledetrois) -- (position);
\draw(regledetrois) -- (pouvoir);
\draw(regledetrois) -- (raison);
\draw(regledetrois) -- (reduire);
\draw(savoir) -- (contraire);
\draw(savoir) -- (demi);
\draw(savoir) -- (donner);
\draw(savoir) -- (falloir);
\draw(savoir) -- (meme);
\draw(savoir) -- (racine);
\draw(science) -- (compter);
\draw(science) -- (dire);
\draw(science) -- (double);
\draw(science) -- (doubler);
\draw(science) -- (ensemble);
%\draw(science) -- (figurer);
\draw(science) -- (garder);
\draw(science) -- (partie);
\draw(science) -- (regle);
\draw(se) -- (chose);
\draw(se) -- (donner);
\draw(se) -- (double);
\draw(se) -- (doubler);
\draw(se) -- (exemple);
\draw(se) -- (garder);
\draw(se) -- (position);
\draw(se) -- (pouvoir);
\draw(se) -- (raison);
\draw(se) -- (regledetrois);
\draw(seulement) -- (chose);
\draw(seulement) -- (exemple);
\draw(seulement) -- (multiplication);
\draw(seulement) -- (pratique);
\draw(seulement) -- (reduire);
\draw(seulement) -- (regle);
\draw(valoir) -- (donner);
\draw(valoir) -- (double);
\draw(valoir) -- (doubler);
\draw(valoir) -- (exemple);
\draw(valoir) -- (position);
\draw(valoir) -- (pouvoir);
\draw(valoir) -- (raison);
\draw(valoir) -- (regledetrois);
\draw(valoir) -- (savoir);
\draw(valoir) -- (se);
\end{tikzpicture}

\caption{Graph of the relations between fickle words. Two nodes are connected if the words are significantly neighbors.}
\label{fig:graph1}
\end{figure}

In the case of fickle words, the graph~$G$ is given by Figure~\ref{fig:graph1}.

\subsection{Quasi-cliques}

Graphs are powerful tools for visualization, since the graphical representation can be built according to some parameters that ensure highly connected set of vertices to be gathered as much as possible. Still, it can be interesting not to rely only on graphical intuition, but also to use some clustering algorithms with performance guarantee.\\

Since our graph is pretty dense, it appears that the concept we need here is a quasi-clique coloring. For an introduction to quasi-clique and clique partition problems, one can for example refer to \cite{AJLR10}. A \textit{quasi-clique} is a subgraph of highest density; typically, if $h$ is a nondecreasing function, $K \subset V$ is a quasi-clique according to $h$ if $|E[K]| \geq h(|K|,|V|)$. Note that we use the following notations, which are classical in graph theory: if $W\subset V$ is a subset of vertices, then $G[W]$ is the subgraph induced by $W$, and $E[W]$ is set of edges which are internal to $G[W]$. Here we choose $h:|K|,|V|\mapsto|K|(|K|-1)/2-1$. In other terms, we define a quasi-clique as a subgraph such that every pair of vertices except at most one is connected.\\

Unfortunately, finding a maximum quasi-clique is NP-hard in the general case as well as with this specific function. Still since the graph is small and has bounded degree, we can afford to use moderately exponential algorithms (see~\cite{qis}). We build a partition of the graph as such:\\

\begin{algorithm}[ht!]
\textbf{\algo{Glutton Quasi-Clique Decomposition}($G$)}
\label{algo:glutton}
\begin{algorithmic}[1]
    \IF {$\kappa = \{ H \subset V$, $|H|\geq 4 \vee E[H] \leq |H|(|H|-1)/2-1) \} \neq \emptyset$}
        \STATE Find $K$ which is maximum among $\kappa$
        \STATE Return $(K,{\rm GLUTTON~QUASI-CLIQUE~DECOMPOSITION}(G[V\setminus K]))$
    \ELSE
    	\STATE Return $V$
    \ENDIF
\end{algorithmic}
\end{algorithm}

This formal definition can be rephrased with a simple explanation: the algorithm will look for the maximum size quasi-clique, add it as an item of the partition and proceed recursively until the graph contains no quasi-clique of size $4$ or more. The remaining vertices are left isolated in the decomposition (note that a quasi-clique of size $3$ is simply a path, and thus not very interesting to study).\\

Of course the difficult point is the computation of $\kappa$ at each step. There are basically two solutions, depending on the density of the graph.\\

If the graph $G(V,E)$ has high average degree $\delta_G = 2|E|/|V|$, then its complementary graph $\bar{G}(V,V^2\setminus E)$ has low average degree $\delta_{\bar{G}} = |V|-1-\delta_G$, and thus the algorithm from~\cite{qis} is efficient.\\

On the other hand, if the graph has low average degree, we can use the following algorithm, that basically solves the quasi-clique problem through the resolution of a small (quadratic) number of clique problems:

\begin{algorithm}[ht!]
\textbf{\algo{QUASI CLIQUE}($V,E$)}
\label{algo:param}
\begin{algorithmic}[1]
	\STATE $K = {\rm CLIQUE}(V,E)$
	\FORALL{$u,v \in V (u,v)\notin E$}
		\STATE {$K = \max\{K,{\rm CLIQUE}(V,E\cup(u,v))\}$}
	\ENDFOR
	\RETURN $K$
\end{algorithmic}
\end{algorithm}

Here CLIQUE can be any exact algorithm for the maximum clique problem, which is \textbf{NP}-hard too. To our knowledge, the fastest ones are those designed in~\cite{exactis}.\\

\section{Analysis of results}\label{analysis}

\subsection{Robust Kohonen maps}

In what follows, we consider the Kohonen map after removing the fickle words (in gray), see Figure~\ref{fig:koho3}. We call this modified map a Robust Kohonen map. 

The Robust Kohonen map shows a contrast between the top right corner and the bottom left one, the same contrast as between left and right sides on the first axis in FCA representation (see Figure~\ref{fig:CA12}).

Indeed, the top right corner contains words linked to arithmetic practice and verbs that are used to build arithmetical operations as \textit{retenir}"to retain", \textit{emprunter}"to borrow", etc. No manuscript is specific of this part of the map, even if the BNF fr. 2050 and \textit{Kadran aux Marchans} appear in the low part of the graph. 

On the bottom left corner, one finds the lexical inheritance of the medieval university  represented by the manuscript BNF fr. 1339. Some texts contain a highly specialized vocabulary, with connections to the university world, with words such as \textit{article} "article", \textit{algorism} "algorithm", \textit{sain} "integer", or a vocabulary of geometry used for the extraction of roots (\textit{carree}).\\

The other two corners of the map are characterized by vocabulary taken from two specific manuscripts. The BSG 3143 in the up left corner is a treatise written by Jean Adam for future Louis XI, that is exceptional in the corpus because it uses Latin words and roman numbers and also because it had to be pleasant for the prince. In spite of this, it shares with the Nantes 456 and BNF fr.1339 words as \textit{gecter}, \textit{gectons} that are marks of abacus algorithms.

In the opposite corner, the \textit{Traicté de la pratique} is marked by a more descriptive vocabulary of mathematical problems (\textit{item} "item", \textit{demande} "demand", \textit{requerir} "to call for", \textit{quant} "quant")  and stronger scientific approach (\textit{aliquot} "aliquot", \textit{corps} "field", \textit{proportionnellement} "proportionally".)\\

What can we do with the list of "fickle words" from this map ? First, it is remarkable that a part of fickle words concerns the algorithm of the rule of three. This algorithm  consists of a "multiplication" (\textit{multiplier}) by the "opposite" (\textit{contraire}) and of a "division" (\textit{diviser}). Other fickle words are related to the operations (\textit{reduction} "fractions reducing", \textit{multiplier} "to multiply", \textit{additionner} "to add"), and with words having a distinctive didactic flavor (\textit{falloir} "have to", \textit{dire} "to say"). 
As a matter of fact, the two main technical issues for these XV\up{th} century authors are to teach how to use the rule of three and the fractions to their readers.

\hyphenpenalty=10000
\setlength{\tabcolsep}{1mm}
\setlength{\arrayrulewidth}{1pt}

\begin{figure}[ht!]
\centering
\scalebox{0.6}{ 
\begin{tabular}{|p{1.9cm}|p{1.9cm}|p{1.9cm}|p{1.9cm}|p{1.9cm}|p{1.9cm}|p{1.9cm}|p{1.9cm}|p{1.9cm}|p{1.9cm}|}
\hline
minutes super* & & notes & calculer cubic & fois mettre \textbf{BNF10259} & &  \textcolor{gris}{contraire}  \textcolor{gris}{depenser}  \textcolor{gris}{falloir}  \textcolor{gris}{meme}  \textcolor{gris}{racine} & aller donc ensuivre  \textcolor{gris}{savoir} &  \textcolor{gris}{multiplier} &  \textcolor{gris}{regle} venir\\
\hline
gecter & & & & & dessous & barrater  \textcolor{gris}{demi} & somme & voir & assembler\\
\hline
 &  \textbf{BSG3143} & parteur & defaillir duplation mediation nommer numeration senestre & semblables & & circulaires demeurer derenier disaine ecrire entendre formes nombrateur oter prouver & entrer laisser rien & emprunter figure regarder & figure de~non~rien fraction muer rayes retenir\\
\hline
notables & & denomi- nations multipli- cateur nominateur ordonne & abaisser comptes endroit proposer & anteriorer diminution enseignement enseigner moyen signifiant surplus trancher & possible  \textcolor{gris}{reduire} & difficile progression repondre & avaluer & faillir & monter partiteur\\
\hline
bref gectons  \textcolor{gris}{multipli-} \textcolor{gris}{cation}  \textcolor{gris}{pratique}  \textcolor{gris}{seulement} & generale latin  \textcolor{gris}{nombrer} proportion reduction unite & soustraction & denomi- nateur & entier & ajoutement & & remotion \textbf{Kadran} & & \textbf{BNF2050}\\
\hline
chose& nulle &  \textcolor{gris}{ensemble}  \textcolor{gris}{partie} & \textbf{Nantes456} & & & &  \textcolor{gris}{partement}  \textcolor{gris}{valoir} & & etre\\
\hline
arithmetique  \textcolor{gris}{compter} preuve tenir & cubbement destre digit diviser diviseur division lignes nombre compo \textcolor{gris}{se}querir &  \textcolor{gris}{dire}  \textcolor{gris}{figurer} poser  \textcolor{gris}{science} & & & abreger lever precedent quotiens & numerateur & moindre nombre \textbf{Nicolas Chuquet} & partir plus  \textcolor{gris}{raison} & \\
\hline
egalir egaliser especes question total traiter & mesurer & grand & apparaitre  \textcolor{gris}{position}  \textcolor{gris}{soustraire} &  \textcolor{gris}{ajouter} & leurs prendre quant quantefois reponse trouver & &  \textcolor{gris}{part} & commun Item & devoir  \textcolor{gris}{droit}  \textcolor{gris}{exemple} reste rester\\
\hline
algorisme article carrees cercle envient ligne pair sain \textbf{BNF1339} & addition chiffre former &  \textcolor{gris}{double}  \textcolor{gris}{doubler} moitie & appeler  \textcolor{gris}{donner}  \textcolor{gris}{maniere}  \textcolor{gris}{pouvoir}  \textcolor{gris}{se}& bailler demander mises nomper pareillement vouloir & & egale faire montrer necessaire romp selon & & & \\
\hline
cautelle & demontrer dessus &  \textcolor{gris}{garder}  \textcolor{gris}{regle}  \textcolor{gris}{de}  \textcolor{gris}{trois} & aliquot composer corps moins sub* toutefois & partant plaisir proportionel- lement requerir residu & appartenir convenir demande difference egaulx maieur millions rate survendre tant & & naturel roupt \textbf{Traicte praticque} & & fausse\\
\hline
\end{tabular}
}
\caption{Robust Kohonen map: fickle words are removed (in gray)}
\label{fig:koho3}
\end{figure}

\subsection{Improved visualization for FCA}

\begin{figure}
\centering
\begin{tikzpicture}[scale=0.14]

\pgfmathsetmacro\fsize{87};
\pgfmathsetmacro\xorg{20};
\pgfmathsetmacro\yorg{20};
\pgfmathsetmacro\xdlcorner{-11};
\pgfmathsetmacro\ydlcorner{-23};
\pgfmathsetmacro\stepunit{25};

\coordinate (dlcorner) at (\xdlcorner,\ydlcorner);

\draw[black,line width=0.4mm] (dlcorner) -- +(0,\fsize) -- +(\fsize,\fsize) -- +(\fsize,0) -- +(0,0);
\draw[dotted,black,line width=0.4mm] (\xdlcorner+2.5,\yorg) -- (\xdlcorner+\fsize,\yorg);
\draw[dotted,black,line width=0.4mm] (\xorg,\ydlcorner+1.8) -- (\xorg,\ydlcorner+\fsize);

\begin{scope}[every node/.style={font={\scriptsize}}]
\node at (32,\ydlcorner+1){\textbf{ Factor 1 (25.03\%)}};
\node at (\xorg-\stepunit,\ydlcorner+1){\textbf{-1}};
\node at (\xorg,\ydlcorner+1){\textbf{0}};
\node at (\xorg+\stepunit,\ydlcorner+1){\textbf{1}};
\node at (\xorg+2*\stepunit,\ydlcorner+1){\textbf{2}};
\node[rotate=90] at (\xdlcorner+1.7,7){\textbf{ Factor 2 (17.91\%)}};
\node at (-9.6,\yorg){\textbf{0}};
\node at (-9.4,\yorg-\stepunit){\textbf{-1}};
\node at (-9.6,\yorg+\stepunit){\textbf{1}};
\end{scope}

\path plot[mark=*,mark size=4mm] coordinates {(18.375,19.375) (28.875,31.425) (27.225,26.475) (15.925,32.55) (31.775,26.05) (34.425,32.2) (27.875,7.675) (22.375,26.6) (30.1,29.4) (29.25,30.075) (15.525,21.15) (33.1,23.375) (16.375,21.4) (32.625,25.625) (25.125,6.175) (31.65,20.6) (27.15,21.375) (30.475,28.6) (13.9,26.05) (24.35,24.975) (24.45,21.725) (36.975,30.725) (18.825,30.35) (15.2,16.9) (30.175,25.275) (24.8,17.875) (24.55,16.8) (20.125,24.175) (41.9,35.1) (16.25,15.775) (27.975,24.075) (25.75,14.475) (24.15,24.975) (18.375,21.5) (31.525,12.725) (19.4,21.925) (21.55,27.775) (25.675,17.65) (20.95,20.25) (17.425,10.775) };

\path plot[mark=o,mark size=4mm] coordinates {(50.625,17.75) (17.775,16.475) (37.675,29.475) (37.125,14.4) (47.575,27.325) (4.225,40.025) (22.55,30.925) (44.05,9.825) (26.575,25.225) (11.15,34.125) (21.5,28.65) (42.025,23.15) (51.925,24.15) (20.875,17.45) (11.525,21.525) (13.1,30) (26.25,16.1) (27.3,26.3) (27.2,46.1) (60.275,45.45) (60.575,40.95) (58.6,40.35) (49.9,26.525) (45.225,23.175) (14.9,21.25) (23.85,33.1) (40.175,11.375) (5.725,39.8) (25.425,35.875) (55.85,16.225) (35.875,50.6) (44.525,15.1) (6.975,33.675) (11.45,33.2) (36.8,30.375) (37.5,40.425) (26.65,12.075) (41.125,18.675) (22.25,18) (32.6,35.975) (41.975,36.925) (48.575,14.075) (13.8,15.725) (4.2,36.925) (26,11.35) (49.475,12.4) (40,6.1) (26.05,20.15) (52.675,16.725) (50.4,26.8) (44.75,18.675) (22.3,37.175) (47.825,19.175) (31.375,24.375) (8.325,27.55) (51.8,26.35) (61.35,26.375) (-0.125,32.925) (15.35,7.8) (31.875,17.7) (40.05,10.775) (40.625,10.125) (12.4,32.95) (30,19.725) (33.675,9.225) (18.15,5.9) (71.925,40.625) (56.075,26.325) (16.1,10.1) (16.125,11.275) (11.35,25.375) (11,12.5) (17,-0.775) (9.25,-20.275) (30.825,45.9) (47.25,36.85) (43.325,25.05) (22.325,-2.175) (38.35,33.9) (8.525,37.9) (32.925,19.35) (35.625,18.575) (17.05,12.975) (25.925,14.325) (41.6,33.8) (7.775,28.025) (20.575,11.1) (55.675,37.225) (49.15,20.425) (-6.9,36.825) (50.75,23.25) (42.6,30.475) (33.55,44.075) (2.35,35.625) (44.55,28.55) (18.525,31.125) (27.15,22.65) (19.1,28.525) (36.95,24.15) (8.725,7.95) (38.225,15.4) (9.725,-12.075) (39.65,18.95) (7.075,16.45) (9.05,26.85) (20.75,15.5) (18.925,21.05) (49.15,16.025) (45.05,38.075) (34.2,16.975) (18.5,37.25) (18.25,39) (37.175,43.875) (39.825,28.275) (21.4,24.175) (46.5,22.525) (28.925,33.325) (34.5,23.55) (56.35,33.9) (19.525,38.025) (6.875,28.775) (42.3,24.6) (15.125,12.45) (15.875,1.95) (1.825,31.425) (15.575,19.625) (26.775,9.85) (30.4,13.325) (28.675,16.5) (17.225,27.025) (29.875,25.425) (27.075,19.3) (36.925,29.25) (2.425,38.725) (33.05,21.4) (26.05,16.525) (12.325,25.95) (12.6,28.475) (49.45,22.7) (41.25,31.475) (23.65,23.175) (-4.625,36.825) (8.1,-22.25) (29.85,20.975) (16.975,6.95) (20.875,23.775) (19.8,22.775) (9.05,27.65) (-0.425,33.4) (10.525,39.875) (8.225,11.95) (7.875,17.125) (23.175,1.65) (24.3,10.675) (10.675,23.5) (5.775,16.6) (64.45,38.275) (12.45,26.65) (30.15,21.975) (31.475,15.975) (50.2,8.175) (28.65,19.175) (31.35,17.375) (12.7,36.35) (41.375,40.6) (34.775,10.6) (-2.55,35.95) (3.7,31.525) (37.4,18.975) (43.375,17.55) (11.425,36.7) (36.75,28.65) (44.875,3.125) (8.225,27.6) (48.975,29.65) (15.275,18.2) (13.95,28.85)};

\begin{scope}[every node/.style={font={\tiny},inner sep=1pt,fill=white}]
\node[draw] (Kadran) at (18.95,15.075){{Kadran}};
\node[draw] (BSG3143) at (28.55,26.925){{BSG3143}};
\node[draw] (BNF1339) at (41.025,26.95){{BNF1339}};
\node[draw] (Nicolas-Chuquet) at (13.675,26.7){{Nicolas Chuquet}};
\node[draw] (BNF2050) at (15.525,5.5){{BNF2050}};
\node[draw] (BNF10259) at (25.825,30.975){{BNF10259}};
\node[draw] (Traicte-praticque) at (9.175,25.45){{Traicte praticque}};
\node[draw] (Nantes456) at (29.825,14.125){{Nantes456}};
\end{scope}

\node (dire)at (27.875,7.675){};
\node (raison)at (16.25,15.775){};
\node (regle-de-trois)at (24.15,24.975){};
\node (pratique)at (20.125,24.175){};

\begin{scope}[every node/.style={font={\footnotesize}}]
\draw[black,<-,line width=0.2mm] (dire) -- +(4,-1)node[right]{dire};
\draw[black,<-,line width=0.2mm] (raison) -- +(-4,-2)node[left]{raison};
\draw[black,<-,line width=0.2mm] (regle-de-trois) -- +(4,-1)node[right]{règle de trois};
\draw[black,<-,line width=0.2mm] (pratique) -- +(-4,-2)node[left]{pratique};
\end{scope}

\end{tikzpicture}
\caption{Projection on first two factors of the FCA. Only the fickle words are in black.}
\label{fig:CA12fi}
\end{figure}
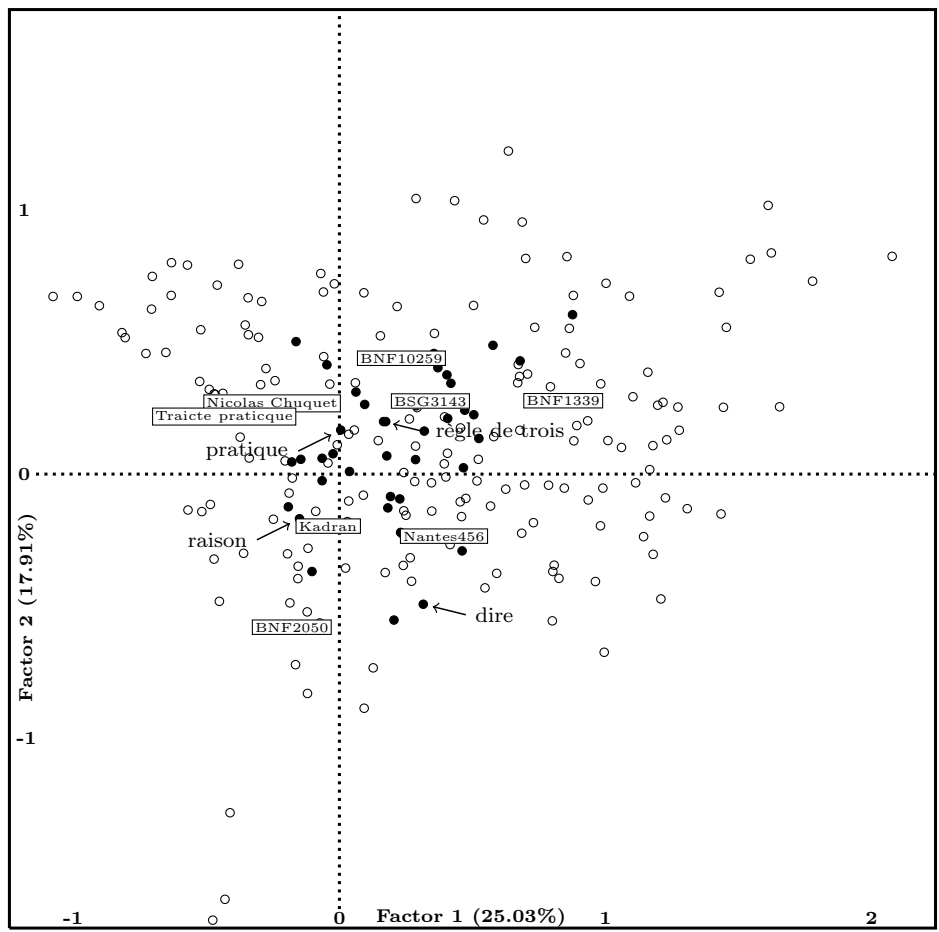

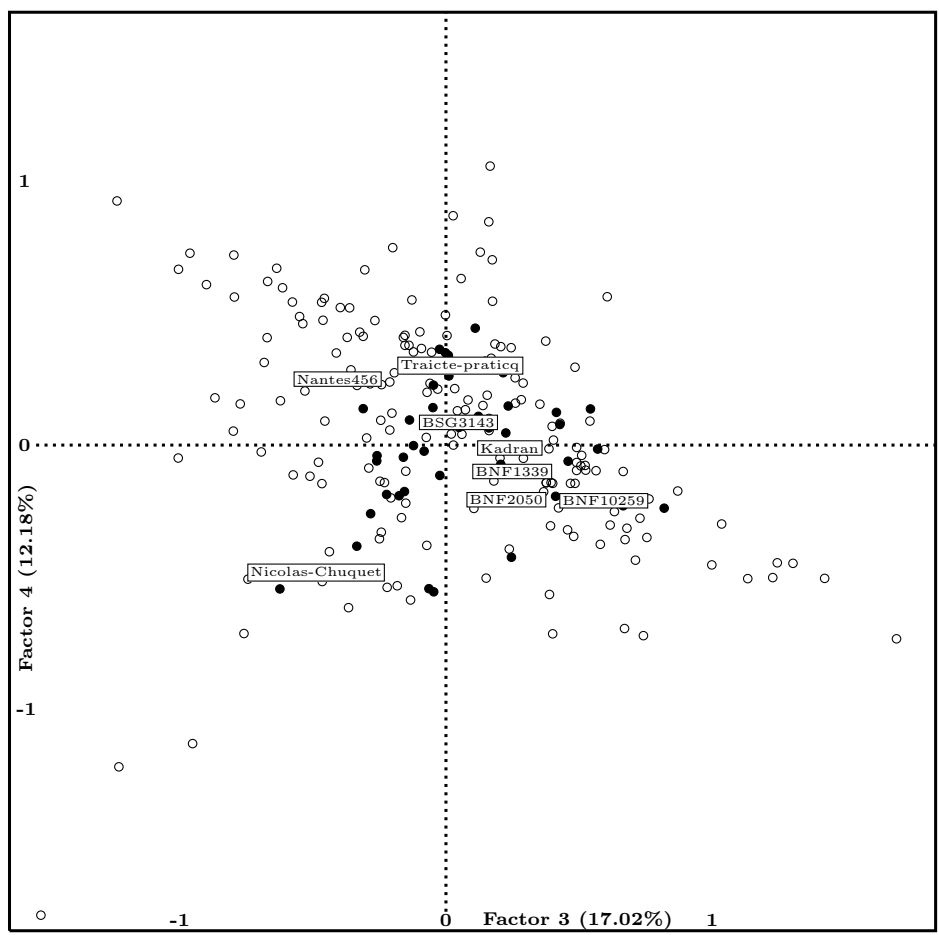
\begin{figure}[h!]
\centering
\begin{tikzpicture}[scale=0.14]

\pgfmathsetmacro\fsize{87};
\pgfmathsetmacro\xorg{20};
\pgfmathsetmacro\yorg{20};
\pgfmathsetmacro\xdlcorner{-21};
\pgfmathsetmacro\ydlcorner{-26};
\pgfmathsetmacro\stepunit{25};

\coordinate (dlcorner) at (\xdlcorner,\ydlcorner);

\draw[black,line width=0.4mm] (dlcorner) -- +(0,\fsize) -- +(\fsize,\fsize) -- +(\fsize,0) -- +(0,0);
\draw[dotted,black,line width=0.4mm] (\xdlcorner+2.5,\yorg) -- (\xdlcorner+\fsize,\yorg);
\draw[dotted,black,line width=0.4mm] (\xorg,\ydlcorner+1.8) -- (\xorg,\ydlcorner+\fsize);

\begin{scope}[every node/.style={font={\scriptsize}}]
\node at (32,\ydlcorner+1){\textbf{ Factor 3 (17.02\%)}};
\node at (\xorg-\stepunit,\ydlcorner+1){\textbf{-1}};
\node at (\xorg,\ydlcorner+1){\textbf{0}};
\node at (\xorg+\stepunit,\ydlcorner+1){\textbf{1}};
\node[rotate=90] at (\xdlcorner+1.7,7){\textbf{ Factor 4 (12.18\%)}};
\node at (-19.6,\yorg){\textbf{0}};
\node at (-19.4,\yorg-\stepunit){\textbf{-1}};
\node at (-19.6,\yorg+\stepunit){\textbf{1}};
\end{scope}

\path plot[mark=o,mark size=4mm] coordinates {(5.575,33.55) (5.65,17.175) (29.675,19.65) (16.95,28.825) (32.275,18.35) (-18.05,-24.525) (25.075,18.775) (-5.125,36.65) (4.45,24.2) (24.35,37.55) (16.225,17.525) (32.675,18.025) (30.025,16.375) (23.775,7.4) (29.375,29.85) (16.225,14.5) (44.975,8.65) (20.7,20) (52.6,8.8) (48.35,7.35) (33.125,17.625) (36.65,17.5) (32.125,16.375) (39.075,14.9) (9.05,9.9) (15.825,13.125) (-2.5,35.2) (23.875,24.725) (-3.8,-8.275) (6.55,31.5) (55.575,7.375) (0.075,38) (24.25,28.225) (24.05,21.375) (32.75,19.025) (15.425,6.675) (16,30.2) (8.325,33.525) (32.275,17.6) (51.125,8.85) (25.95,10.15) (8.45,31.825) (34.9,19.575) (35.15,34.05) (27.25,25.875) (9.7,28.725) (11.075,27.125) (29.175,15.575) (6.25,32.175) (21.825,23.35) (12.55,20.675) (33.525,22.275) (12.225,30.3) (34.1,17.575) (16.8,33.75) (18.15,20.725) (18.575,22.225) (15,38.7) (29.725,5.85) (3.2,30.175) (4.65,34.9) (3.25,35.5) (26.525,23.975) (27.075,24.3) (14.725,25.975) (38.225,13.075) (38.875,11.25) (19.225,25.3) (20.8,27.9) (28.825,23.875) (22.075,24.275) (24.6,29.575) (36.825,11.05) (38.55,1.95) (62.325,1.65) (30.725,22.075) (30.1,20.475) (29.825,12.35) (20.1,30.375) (24.025,41.15) (20.15,27.35) (13.8,16.575) (16.675,5.325) (30.575,14.05) (32.275,19.775) (8.625,22.275) (8.375,7.075) (29.425,16.45) (6.75,25.125) (20.675,41.725) (11.9,30.7) (32,11.35) (50.7,7.45) (19.95,32.325) (13.325,31.8) (21.5,21.025) (2.65,19.35) (14.475,6.525) (14.225,16.45) (33.075,18.075) (4.1,36.75) (31.425,11.975) (8.575,33.9) (21.425,35.775) (20.825,25.35) (45.9,12.525) (8.35,16.35) (10.725,30.2) (29.975,21.775) (16.15,29.45) (21.075,23.25) (18.65,28.8) (41.775,15.65) (15.15,26.85) (-1.7,24.475) (10.075,33.025) (20.5,21.05) (29.425,16.425) (24.475,19.975) (23.475,23.75) (13.875,22.35) (10.95,33) (10.85,4.6) (29.875,16.45) (8.025,18.375) (13.925,11.75) (17.7,29.15) (22.35,27.25) (0.675,23.9) (7.225,17.05) (26.5,26.375) (25.15,29.325) (17.55,30.725) (1.4,7.3) (2.925,27.825) (22.625,14) (14.925,23.025) (0.025,21.325) (21.8,22.45) (31.7,16.375) (-5.15,18.775) (24.15,46.425) (36.775,2.625) (12.375,36.6) (30.025,2.125) (24.375,33.625) (32.125,27.375) (13.75,11.125) (12.875,25.8) (14.825,15) (27.275,18.75) (18.2,10.5) (34.5,10.6) (37,12.125) (16.525,29.45) (26.125,29.225) (35.425,12.425) (18.5,25.85) (11.65,25.65) (13.95,25.725) (-4.05,38.175) (29.375,17.275) (16.15,30.4) (-10.725,-10.475) (23.675,27.95) (0.125,34.025) (23.225,38.275) (18.225,25) (24.5,16.6) (37.8,9.1) (1.025,2.15) (35.825,13.7) (-10.9,43.125) (12.75,17.825) (14.725,21.425) (24.05,22.525) (20.275,22.05)};
\path plot[mark=*,mark size=4mm] coordinates { (15.6,15.2) (25.15,18.175) (22.75,31.075) (36.65,14.25) (34.25,19.625) (30.3,15.15) (17.95,19.425) (19.425,17.125) (13.525,19) (12.225,23.45) (33.575,23.425) (19.95,28.75) (23.05,22.725) (26.15,9.375) (16,18.85) (12.925,13.5) (4.4,6.375) (31.475,18.475) (20.275,26.55) (30.7,21.95) (25.625,21.15) (30.375,23.1) (18.825,25.675) (16.575,22.375) (19.4,29.075) (14.425,15.325) (18.4,6.4) (20.225,28.5) (40.5,14.025) (18.85,6.1) (25.375,26.85) (25.85,23.7) (11.625,10.425) (24.025,21.625) (21.225,21.675) (16.1,15.6) (18.775,23.55) (13.5,18.5) (16.95,19.95) (28.3,19.725)};

\begin{scope}[every node/.style={font={\tiny},inner sep=1pt,fill=white}]
\node[draw] (Kadran)at (26,19.725){{Kadran}};
\node[draw] (BSG3143)at (21.175,22.15){{BSG3143}};
\node[draw] (BNF1339)at (26.225,17.5){{BNF1339}};
\node[draw] (Nicolas-Chuquet)at (7.8,7.925){{Nicolas-Chuquet}};
\node[draw] (BNF25)at (25.675,14.825){{BNF2050}};
\node[draw] (BNF1259)at (34.825,14.725){{BNF10259}};
\node[draw] (Traicte-praticq)at (21.35,27.525){{Traicte-praticq}};
\node[draw] (Nantes456)at (9.825,26.25){{Nantes456}};
\end{scope}

\end{tikzpicture}
\caption{Projection on third and fourth factors of the FCA. Only the fickle words are in black.}
\label{fig:CA34fi}
\end{figure}

The combination of both techniques FCA and SOM whose result is displayed in Figures~\ref{fig:CA12fi} and~\ref{fig:CA34fi} is interesting because it preserves properties of the FCA while giving additional information about the center of the projection - which is usually difficult to interpret. Indeed, the identification  of the fickle words on the FCA projections allows us to improve the general interpretation of the factorial graphs, where some words are located because of the algorithm and not because of their attraction to other words and to the texts.\\

Remember that, on the first two factorial axes (see Figure~\ref{fig:CA12fi}), we have observed an opposition between the university legacy, on the right, and a more practical pole with rules, problems and fractions, on the left. It could be tempting to support this observation with very significant words such as \textit{pratique} "practical" or \textit{règle de trois} "rule of three". Still, the enhancement of the fickle forms on the FCA shows that these words are in fact shared between many different texts and not only linked to the more 'practical' ones: Nicolas Chuquet and \textit{Traict\'e en la praticque}. As a matter of fact, they do belong to all the texts.\\ 

It is the same for two other words (\textit{raison} "reason", \textit{dire} "to say"). The word \textit{raison} is an ambiguous word, in a way, because it can mean calculation, with textual matches like "do your reasons", or indicate mathematical problems. "To say" ranks sixth by order of frequency among verbs in the corpus. 
Note that all most important verbs are not fickle words. The first eight verbs for occurrence are: \textit{être}(14523) "to be", \textit{avoir}(4563) "to have", \textit{devoir}(3826) "must", \textit{faire}(3431) "to do", \textit{multiplier}(3228) "to multiply", \textit{dire}(2461) "to say", \textit{partir}(2648) "to divide", "valoir"(2606) "to be worth". One of the particular meanings of "to say" comes from the orality of this type of text. Understanding arithmetic operations often supposes saying it aloud.\\

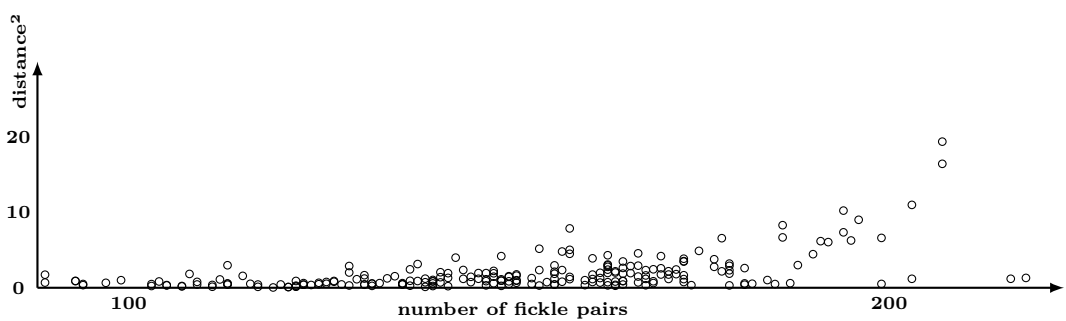
\begin{figure}[h!]
\centering
\begin{tikzpicture}[scale=0.1]
\coordinate (origin) at (88,0);
\draw[thick,black,-latex](origin)-- +(135,0);
\draw[thick,black,-latex](origin)-- +(0,30);
\path plot[mark=o,mark size=5mm] coordinates {(89,1.733313)
(89,0.687671)
(93,0.879188)
(93,0.917328)
(94,0.500639)
(94,0.349571)
(97,0.657271)
(99,1.005228)
(103,0.488729)
(103,0.233705)
(104,0.834061)
(105,0.355916)
(105,0.311233)
(107,0.179992)
(107,0.236564)
(108,1.834052)
(109,0.780884)
(109,0.433761)
(111,0.374555)
(111,0.136533)
(112,1.102961)
(113,0.549119)
(113,0.563886)
(113,0.407632)
(113,2.986735)
(115,1.576985)
(116,0.525266)
(117,0.12585)
(117,0.427956)
(119,0.039282)
(120,0.43786)
(121,0.089967)
(121,0.106757)
(122,0.904057)
(122,0.27663)
(122,0.143855)
(122,0.170589)
(123,0.408344)
(123,0.409868)
(123,0.619592)
(124,0.370105)
(125,0.662793)
(125,0.470098)
(126,0.53128)
(126,0.771339)
(127,0.902703)
(127,0.799271)
(128,0.443143)
(129,0.302426)
(129,2.023817)
(129,2.884816)
(130,1.131654)
(131,1.667429)
(131,0.454524)
(131,0.402024)
(131,1.22218)
(132,0.579035)
(132,0.526653)
(132,0.251264)
(133,0.616724)
(134,1.246114)
(135,1.529439)
(136,0.61385)
(136,0.49528)
(136,0.433612)
(136,0.45828)
(137,0.273524)
(137,0.925961)
(137,2.444701)
(138,3.146146)
(138,0.895174)
(139,1.207906)
(139,0.730987)
(139,0.122534)
(140,0.86976)
(140,0.582069)
(140,0.833297)
(140,1.039637)
(140,0.230508)
(141,1.403154)
(141,2.060198)
(141,0.709805)
(142,1.934577)
(142,0.211917)
(142,1.309469)
(143,3.995438)
(144,1.184564)
(144,2.358767)
(145,0.747727)
(145,1.454136)
(146,1.193686)
(146,1.94952)
(147,1.944066)
(147,1.130539)
(147,0.318043)
(147,0.879615)
(148,1.468539)
(148,2.237699)
(148,1.919718)
(148,0.593212)
(149,1.119967)
(149,4.19534)
(149,0.943449)
(149,0.253371)
(149,0.939117)
(150,1.532861)
(150,1.440959)
(150,0.496267)
(150,0.853435)
(151,0.668082)
(151,1.773085)
(151,0.934419)
(151,1.001479)
(151,1.609046)
(153,0.502958)
(153,1.281498)
(154,0.294367)
(154,5.175975)
(154,2.346217)
(155,0.746928)
(156,0.412462)
(156,1.195406)
(156,2.969508)
(156,1.857245)
(156,2.124791)
(156,0.44429)
(157,2.346391)
(157,0.801633)
(157,4.806805)
(158,1.442624)
(158,4.504724)
(158,7.870092)
(158,1.14631)
(158,5.032216)
(160,0.350625)
(160,1.006146)
(161,0.770879)
(161,3.903592)
(161,1.154718)
(161,1.745997)
(162,1.985234)
(162,0.689118)
(162,1.289477)
(163,2.403167)
(163,0.320264)
(163,4.312527)
(163,2.970866)
(163,3.079829)
(163,1.966399)
(163,2.773674)
(163,0.741294)
(164,2.02266)
(164,2.180347)
(164,0.384743)
(164,0.824953)
(164,0.268831)
(165,2.634443)
(165,0.428705)
(165,1.877663)
(165,0.910697)
(165,0.873904)
(165,3.485168)
(166,2.878361)
(166,1.943489)
(167,2.930207)
(167,1.489835)
(167,0.713504)
(167,4.566103)
(168,1.635471)
(168,2.339072)
(168,1.176074)
(168,0.33371)
(169,2.449229)
(169,0.623658)
(169,0.86062)
(170,1.760697)
(170,4.197148)
(170,2.58299)
(171,1.168305)
(171,2.191193)
(171,1.754257)
(172,1.807177)
(172,2.422776)
(173,1.15045)
(173,3.85446)
(173,0.76326)
(173,1.661414)
(173,3.472332)
(173,3.53184)
(174,0.345418)
(175,4.891938)
(177,3.78168)
(177,2.77811)
(178,6.57394)
(178,2.177587)
(179,2.90853)
(179,0.315672)
(179,2.319238)
(179,3.205815)
(179,1.872733)
(181,0.579635)
(181,2.608949)
(181,0.407873)
(182,0.537897)
(184,1.024579)
(185,0.493946)
(186,6.677435)
(186,8.299131)
(187,0.606432)
(188,3.003859)
(190,4.45374)
(191,6.180186)
(192,6.054028)
(194,7.350671)
(194,10.232712)
(195,6.269358)
(196,9.013623)
(199,0.50955)
(199,6.61164)
(203,1.191249)
(203,10.989775)
(207,16.435334)
(207,19.376364)
(216,1.193709)
(218,1.307591)};

\begin{scope}[every node/.style={font={\scriptsize}}]
\node at (150,-3){\textbf{ number of fickle pairs}};
\node at (100,-2){\textbf{100}};
\node at (200,-2){\textbf{200}};
\node[rotate=90] at (85.5,30){\textbf{distance\up{2}}};
\node at (85.5,0){\textbf{0}};
\node at (85.5,10){\textbf{10}};
\node at (85.5,20){\textbf{20}};
\end{scope}

\end{tikzpicture}
\caption{Correlation between fickleness and distance to the center. x-axis represents the number of fickle pairs a word belongs to, while y-axis stands for the square distance to the origin.}
\label{distances-fickles}
\end{figure}

Another interest of this kind of representation is the interpretation of the center of FCA. As we can see in Figure~\ref{distances-fickles}, fickle words are close to the center, but there are other words in the same place, which we could interpret.\\

To conclude, we can see that two levels of interpretation are superimposed: the fickle pairs reveal the shared lexicon and the factorial map inserts them in a local interaction system. And since the fickle words list is computed independently from the FCA, we can successively study these interactions on each axis. It is the articulation between these two levels which makes this representation interesting. At the end, the meaning of this new kind of factorial map is quite intuitive and offers easy tools to the argumentation.

\begin{table}[h!]
\centering
\setlength{\tabcolsep}{1pt}
\setlength{\arrayrulewidth}{0.5pt}
\scalebox{0.72}{ 
\begin{tabular}{|p{2pt}l|*{30}{c|}}
\hline
& & \begin{sideways}contraire\end{sideways} & \begin{sideways}depenser\end{sideways} & \begin{sideways}falloir\end{sideways} & \begin{sideways}racine\end{sideways} & \begin{sideways}meme\end{sideways} & \begin{sideways}demi\end{sideways} & \begin{sideways}savoir\end{sideways} & \begin{sideways}valoir\end{sideways} & \begin{sideways}regledetrois\end{sideways}&\begin{sideways}raison\end{sideways}&\begin{sideways}position\end{sideways}&\begin{sideways}pouvoir\end{sideways}&\begin{sideways}garder\end{sideways}&\begin{sideways}donner\end{sideways}&\begin{sideways}se\end{sideways}&\begin{sideways}doubler\end{sideways}&\begin{sideways}double\end{sideways}&\begin{sideways}partie\end{sideways}&\begin{sideways}science\end{sideways}&\begin{sideways}compter\end{sideways}&\begin{sideways}ensemble\end{sideways}&\begin{sideways}dire\end{sideways}&\begin{sideways}regle\end{sideways}&\begin{sideways}chose\end{sideways}&\begin{sideways}reduire\end{sideways} & \begin{sideways}seulement\end{sideways} & \begin{sideways}pratique\end{sideways} & \begin{sideways}multiplication \end{sideways}&\begin{sideways}exemple\end{sideways}&\begin{sideways}nombrer\end{sideways} \\ \hline
&contraire&\cellcolor{gray}1&\cellcolor{gray0.55}.5{\tiny 5}&\cellcolor{gray0.5}.6{\tiny 25}&\cellcolor{gray0.5}.6{\tiny 75}&\cellcolor{gray0.45}.7{\tiny 5}&\cellcolor{gray0.5}.6{\tiny 5}&\cellcolor{gray0.6}.4{\tiny 25}&\cellcolor{white}.0{\tiny 75}&\cellcolor{white}.0{\tiny 5}&\cellcolor{white}.0{\tiny 5}&\cellcolor{white}.0{\tiny 25}&\cellcolor{white}.0{\tiny 5}&0&\cellcolor{white}.0{\tiny 75}&\cellcolor{white}.0{\tiny 25}&0&0&\cellcolor{white}.0{\tiny 5}&0&0&\cellcolor{white}.0{\tiny 75}&0&\cellcolor{gray0.9}.1{\tiny 25}&\cellcolor{gray0.9}.1&\cellcolor{gray0.8}.2&\cellcolor{white}.0{\tiny 25}&0&\cellcolor{white}.0{\tiny 75}&\cellcolor{white}.0{\tiny 5}&\cellcolor{white}.0{\tiny 25} \\ \hline
&depenser&\cellcolor{gray0.55}.5{\tiny 5}&\cellcolor{gray}1&\cellcolor{gray0.4}.8{\tiny 5}&\cellcolor{gray0.4}.8&\cellcolor{gray0.45}.7{\tiny 5}&\cellcolor{gray0.55}.5{\tiny 75}&\cellcolor{gray0.9}.1{\tiny 75}&0&\cellcolor{white}.0{\tiny 25}&\cellcolor{white}.0{\tiny 75}&\cellcolor{white}.0{\tiny 5}&\cellcolor{white}.0{\tiny 75}&\cellcolor{white}.0{\tiny 75}&\cellcolor{gray0.9}.1{\tiny 75}&0&\cellcolor{white}.0{\tiny 75}&\cellcolor{gray0.9}.1&\cellcolor{gray0.9}.1{\tiny 5}&\cellcolor{white}.0{\tiny 25}&\cellcolor{gray0.9}.1{\tiny 75}&\cellcolor{gray0.9}.1&\cellcolor{white}.0{\tiny 25}&\cellcolor{white}.0{\tiny 75}&\cellcolor{gray0.9}.1{\tiny 75}&\cellcolor{gray0.9}.1{\tiny 5}&\cellcolor{white}.0{\tiny 25}&\cellcolor{white}.0{\tiny 5}&0&\cellcolor{white}.0{\tiny 5}&\cellcolor{white}.0{\tiny 25} \\ \hline
&falloir&\cellcolor{gray0.5}.6{\tiny 25}&\cellcolor{gray0.4}.8{\tiny 5}&\cellcolor{gray}1&\cellcolor{gray0.5}.6{\tiny 75}&\cellcolor{gray0.45}.7{\tiny 75}&\cellcolor{gray0.55}.5{\tiny 5}&\cellcolor{gray0.8}.2{\tiny 5}&\cellcolor{white}.0{\tiny 5}&\cellcolor{white}.0{\tiny 5}&\cellcolor{gray0.9}.1{\tiny 5}&\cellcolor{gray0.9}.1&\cellcolor{gray0.9}.1&\cellcolor{gray0.9}.1&\cellcolor{gray0.9}.1{\tiny 75}&\cellcolor{white}.0{\tiny 25}&\cellcolor{white}.0{\tiny 5}&\cellcolor{white}.0{\tiny 5}&\cellcolor{gray0.9}.1{\tiny 25}&\cellcolor{white}.0{\tiny 5}&\cellcolor{white}.0{\tiny 75}&\cellcolor{gray0.9}.1{\tiny 5}&\cellcolor{gray0.9}.1{\tiny 25}&\cellcolor{gray0.9}.1&\cellcolor{gray0.9}.1&\cellcolor{gray0.8}.2&0&\cellcolor{white}.0{\tiny 25}&0&\cellcolor{white}.0{\tiny 5}&\cellcolor{white}.0{\tiny 25} \\ \hline
&racine&\cellcolor{gray0.5}.6{\tiny 75}&\cellcolor{gray0.4}.8&\cellcolor{gray0.5}.6{\tiny 75}&\cellcolor{gray}1&\cellcolor{gray0.5}.6{\tiny 75}&\cellcolor{gray0.5}.6{\tiny 25}&\cellcolor{gray0.8}.2{\tiny 5}&\cellcolor{white}.0{\tiny 25}&\cellcolor{white}.0{\tiny 25}&0&0&\cellcolor{white}.0{\tiny 25}&\cellcolor{white}.0{\tiny 25}&\cellcolor{white}.0{\tiny 25}&0&\cellcolor{white}.0{\tiny 25}&\cellcolor{white}.0{\tiny 25}&\cellcolor{white}.0{\tiny 25}&0&0&\cellcolor{white}.0{\tiny 5}&0&\cellcolor{white}.0{\tiny 5}&\cellcolor{white}.0{\tiny 75}&\cellcolor{gray0.9}.1{\tiny 5}&0&0&0&\cellcolor{white}.0{\tiny 25}&\cellcolor{white}.0{\tiny 75} \\ \hline
&meme&\cellcolor{gray0.45}.7{\tiny 5}&\cellcolor{gray0.45}.7{\tiny 5}&\cellcolor{gray0.45}.7{\tiny 75}&\cellcolor{gray0.5}.6{\tiny 75}&\cellcolor{gray}1&\cellcolor{gray0.45}.7{\tiny 25}&\cellcolor{gray0.8}.2{\tiny 75}&\cellcolor{white}.0{\tiny 5}&\cellcolor{white}.0{\tiny 5}&\cellcolor{gray0.9}.1{\tiny 25}&\cellcolor{white}.0{\tiny 75}&\cellcolor{white}.0{\tiny 75}&\cellcolor{white}.0{\tiny 5}&\cellcolor{gray0.9}.1{\tiny 5}&\cellcolor{white}.0{\tiny 25}&\cellcolor{white}.0{\tiny 25}&\cellcolor{white}.0{\tiny 25}&\cellcolor{gray0.9}.1{\tiny 25}&\cellcolor{white}.0{\tiny 25}&\cellcolor{white}.0{\tiny 5}&\cellcolor{gray0.9}.1{\tiny 5}&\cellcolor{white}.0{\tiny 75}&\cellcolor{gray0.9}.1&\cellcolor{gray0.9}.1{\tiny 75}&\cellcolor{gray0.8}.2{\tiny 75}&\cellcolor{white}.0{\tiny 25}&0&\cellcolor{white}.0{\tiny 5}&\cellcolor{white}.0{\tiny 5}&\cellcolor{white}.0{\tiny 25} \\ \hline
&demi&\cellcolor{gray0.5}.6{\tiny 5}&\cellcolor{gray0.55}.5{\tiny 75}&\cellcolor{gray0.55}.5{\tiny 5}&\cellcolor{gray0.5}.6{\tiny 25}&\cellcolor{gray0.45}.7{\tiny 25}&\cellcolor{gray}1&\cellcolor{gray0.7}.3&\cellcolor{white}.0{\tiny 5}&\cellcolor{white}.0{\tiny 5}&\cellcolor{gray0.9}.1{\tiny 75}&\cellcolor{white}.0{\tiny 75}&\cellcolor{white}.0{\tiny 75}&\cellcolor{white}.0{\tiny 75}&\cellcolor{gray0.9}.1{\tiny 5}&\cellcolor{white}.0{\tiny 5}&\cellcolor{white}.0{\tiny 25}&\cellcolor{white}.0{\tiny 25}&\cellcolor{white}.0{\tiny 5}&0&0&\cellcolor{white}.0{\tiny 75}&\cellcolor{white}.0{\tiny 25}&\cellcolor{white}.0{\tiny 5}&\cellcolor{white}.0{\tiny 75}&\cellcolor{gray0.7}.3{\tiny 25}&0&0&\cellcolor{white}.0{\tiny 5}&\cellcolor{white}.0{\tiny 75}&\cellcolor{white}.0{\tiny 25} \\ \hline
&savoir&\cellcolor{gray0.6}.4{\tiny 25}&\cellcolor{gray0.9}.1{\tiny 75}&\cellcolor{gray0.8}.2{\tiny 5}&\cellcolor{gray0.8}.2{\tiny 5}&\cellcolor{gray0.8}.2{\tiny 75}&\cellcolor{gray0.7}.3&\cellcolor{gray}1&\cellcolor{gray0.8}.2&0&\cellcolor{gray0.9}.1&\cellcolor{white}.0{\tiny 75}&\cellcolor{gray0.9}.1&0&\cellcolor{gray0.8}.2&\cellcolor{white}.0{\tiny 5}&\cellcolor{white}.0{\tiny 5}&\cellcolor{white}.0{\tiny 5}&\cellcolor{white}.0{\tiny 75}&\cellcolor{white}.0{\tiny 25}&\cellcolor{gray0.9}.1&0&\cellcolor{white}.0{\tiny 25}&\cellcolor{gray0.9}.1{\tiny 75}&\cellcolor{white}.0{\tiny 75}&\cellcolor{white}.0{\tiny 75}&\cellcolor{gray0.9}.1{\tiny 25}&\cellcolor{gray0.9}.1&\cellcolor{white}.0{\tiny 5}&\cellcolor{gray0.9}.1{\tiny 25}&0 \\ \hline
&valoir&\cellcolor{white}.0{\tiny 75}&0&\cellcolor{white}.0{\tiny 5}&\cellcolor{white}.0{\tiny 25}&\cellcolor{white}.0{\tiny 5}&\cellcolor{white}.0{\tiny 5}&\cellcolor{gray0.8}.2&\cellcolor{gray}1&\cellcolor{gray0.7}.3&\cellcolor{gray0.6}.4{\tiny 25}&\cellcolor{gray0.55}.5{\tiny 25}&\cellcolor{gray0.6}.4{\tiny 5}&\cellcolor{gray0.9}.1{\tiny 25}&\cellcolor{gray0.55}.5{\tiny 5}&\cellcolor{gray0.55}.5{\tiny 5}&\cellcolor{gray0.8}.2&\cellcolor{gray0.8}.2&\cellcolor{gray0.9}.1{\tiny 25}&\cellcolor{white}.0{\tiny 75}&\cellcolor{gray0.9}.1{\tiny 25}&\cellcolor{gray0.9}.1{\tiny 25}&\cellcolor{gray0.9}.1{\tiny 25}&\cellcolor{gray0.9}.1&\cellcolor{white}.0{\tiny 5}&\cellcolor{white}.0{\tiny 5}&\cellcolor{white}.0{\tiny 25}&\cellcolor{gray0.9}.1&0&\cellcolor{gray0.7}.3{\tiny 25}&0 \\ \hline
&regledetrois&\cellcolor{white}.0{\tiny 5}&\cellcolor{white}.0{\tiny 25}&\cellcolor{white}.0{\tiny 5}&\cellcolor{white}.0{\tiny 25}&\cellcolor{white}.0{\tiny 5}&\cellcolor{white}.0{\tiny 5}&0&\cellcolor{gray0.7}.3&\cellcolor{gray}1&\cellcolor{gray0.5}.6{\tiny 75}&\cellcolor{gray0.45}.7&\cellcolor{gray0.55}.5{\tiny 5}&\cellcolor{gray0.5}.6{\tiny 5}&\cellcolor{gray0.6}.4{\tiny 75}&\cellcolor{gray0.5}.6{\tiny 25}&\cellcolor{white}.0{\tiny 75}&\cellcolor{gray0.9}.1{\tiny 75}&\cellcolor{white}.0{\tiny 25}&\cellcolor{white}.0{\tiny 25}&\cellcolor{white}.0{\tiny 25}&0&\cellcolor{white}.0{\tiny 5}&\cellcolor{white}.0{\tiny 25}&\cellcolor{gray0.8}.2{\tiny 25}&\cellcolor{gray0.8}.2{\tiny 25}&\cellcolor{gray0.9}.1{\tiny 25}&\cellcolor{gray0.9}.1&\cellcolor{gray0.9}.1&\cellcolor{gray0.9}.1{\tiny 5}&\cellcolor{white}.0{\tiny 75} \\ \hline
&raison&\cellcolor{white}.0{\tiny 5}&\cellcolor{white}.0{\tiny 75}&\cellcolor{gray0.9}.1{\tiny 5}&0&\cellcolor{gray0.9}.1{\tiny 25}&\cellcolor{gray0.9}.1{\tiny 75}&\cellcolor{gray0.9}.1&\cellcolor{gray0.6}.4{\tiny 25}&\cellcolor{gray0.5}.6{\tiny 75}&\cellcolor{gray}1&\cellcolor{gray0.45}.7{\tiny 25}&\cellcolor{gray0.4}.8&\cellcolor{gray0.6}.4{\tiny 75}&\cellcolor{gray0.7}.3{\tiny 75}&\cellcolor{gray0.5}.6{\tiny 5}&\cellcolor{white}.0{\tiny 25}&\cellcolor{white}.0{\tiny 5}&0&\cellcolor{white}.0{\tiny 25}&0&0&\cellcolor{white}.0{\tiny 5}&0&\cellcolor{white}.0{\tiny 25}&\cellcolor{white}.0{\tiny 5}&\cellcolor{white}.0{\tiny 25}&\cellcolor{white}.0{\tiny 25}&\cellcolor{white}.0{\tiny 25}&\cellcolor{gray0.9}.1{\tiny 5}&0 \\ \hline
&position&\cellcolor{white}.0{\tiny 25}&\cellcolor{white}.0{\tiny 5}&\cellcolor{gray0.9}.1&0&\cellcolor{white}.0{\tiny 75}&\cellcolor{white}.0{\tiny 75}&\cellcolor{white}.0{\tiny 75}&\cellcolor{gray0.55}.5{\tiny 25}&\cellcolor{gray0.45}.7&\cellcolor{gray0.45}.7{\tiny 25}&\cellcolor{gray}1&\cellcolor{gray0.45}.7{\tiny 75}&\cellcolor{gray0.5}.6&\cellcolor{gray0.55}.5{\tiny 25}&\cellcolor{gray0.45}.7{\tiny 25}&\cellcolor{gray0.9}.1{\tiny 75}&\cellcolor{gray0.9}.1{\tiny 75}&\cellcolor{gray0.9}.1&\cellcolor{gray0.9}.1{\tiny 25}&\cellcolor{white}.0{\tiny 5}&\cellcolor{white}.0{\tiny 75}&\cellcolor{gray0.8}.2&0&\cellcolor{gray0.9}.1&\cellcolor{gray0.9}.1&\cellcolor{white}.0{\tiny 5}&\cellcolor{white}.0{\tiny 5}&0&\cellcolor{white}.0{\tiny 75}&0 \\ \hline
&pouvoir&\cellcolor{white}.0{\tiny 5}&\cellcolor{white}.0{\tiny 75}&\cellcolor{gray0.9}.1&\cellcolor{white}.0{\tiny 25}&\cellcolor{white}.0{\tiny 75}&\cellcolor{white}.0{\tiny 75}&\cellcolor{gray0.9}.1&\cellcolor{gray0.6}.4{\tiny 5}&\cellcolor{gray0.55}.5{\tiny 5}&\cellcolor{gray0.4}.8&\cellcolor{gray0.45}.7{\tiny 75}&\cellcolor{gray}1&\cellcolor{gray0.6}.4{\tiny 25}&\cellcolor{gray0.7}.3{\tiny 25}&\cellcolor{gray0.5}.6{\tiny 5}&\cellcolor{white}.0{\tiny 75}&\cellcolor{gray0.9}.1&\cellcolor{white}.0{\tiny 5}&\cellcolor{gray0.9}.1{\tiny 25}&\cellcolor{white}.0{\tiny 25}&\cellcolor{white}.0{\tiny 5}&\cellcolor{gray0.8}.2{\tiny 25}&\cellcolor{white}.0{\tiny 5}&\cellcolor{white}.0{\tiny 25}&0&\cellcolor{white}.0{\tiny 5}&\cellcolor{white}.0{\tiny 5}&\cellcolor{white}.0{\tiny 25}&\cellcolor{white}.0{\tiny 75}&\cellcolor{white}.0{\tiny 25} \\ \hline
&garder&0&\cellcolor{white}.0{\tiny 75}&\cellcolor{gray0.9}.1&\cellcolor{white}.0{\tiny 25}&\cellcolor{white}.0{\tiny 5}&\cellcolor{white}.0{\tiny 75}&0&\cellcolor{gray0.9}.1{\tiny 25}&\cellcolor{gray0.5}.6{\tiny 5}&\cellcolor{gray0.6}.4{\tiny 75}&\cellcolor{gray0.5}.6&\cellcolor{gray0.6}.4{\tiny 25}&\cellcolor{gray}1&\cellcolor{gray0.7}.3{\tiny 25}&\cellcolor{gray0.55}.5&\cellcolor{gray0.8}.2{\tiny 25}&\cellcolor{gray0.7}.3{\tiny 5}&\cellcolor{gray0.8}.2&\cellcolor{gray0.8}.2&\cellcolor{gray0.9}.1{\tiny 75}&\cellcolor{gray0.9}.1{\tiny 75}&\cellcolor{gray0.9}.1&\cellcolor{white}.0{\tiny 5}&\cellcolor{gray0.8}.2&\cellcolor{gray0.7}.3{\tiny 25}&\cellcolor{white}.0{\tiny 5}&\cellcolor{white}.0{\tiny 75}&\cellcolor{white}.0{\tiny 25}&\cellcolor{gray0.9}.1&\cellcolor{white}.0{\tiny 75} \\ \hline
&donner&\cellcolor{white}.0{\tiny 75}&\cellcolor{gray0.9}.1{\tiny 75}&\cellcolor{gray0.9}.1{\tiny 75}&\cellcolor{white}.0{\tiny 25}&\cellcolor{gray0.9}.1{\tiny 5}&\cellcolor{gray0.9}.1{\tiny 5}&\cellcolor{gray0.8}.2&\cellcolor{gray0.55}.5{\tiny 5}&\cellcolor{gray0.6}.4{\tiny 75}&\cellcolor{gray0.7}.3{\tiny 75}&\cellcolor{gray0.55}.5{\tiny 25}&\cellcolor{gray0.7}.3{\tiny 25}&\cellcolor{gray0.7}.3{\tiny 25}&\cellcolor{gray}1&\cellcolor{gray0.5}.6{\tiny 75}&\cellcolor{gray0.6}.4{\tiny 5}&\cellcolor{gray0.6}.4{\tiny 75}&\cellcolor{gray0.7}.3&\cellcolor{white}.0{\tiny 75}&\cellcolor{gray0.7}.3{\tiny 75}&\cellcolor{gray0.9}.1{\tiny 75}&\cellcolor{white}.0{\tiny 25}&\cellcolor{white}.0{\tiny 25}&\cellcolor{gray0.9}.1{\tiny 5}&\cellcolor{white}.0{\tiny 75}&\cellcolor{gray0.9}.1{\tiny 75}&\cellcolor{gray0.9}.1{\tiny 25}&\cellcolor{white}.0{\tiny 25}&\cellcolor{gray0.7}.3&0 \\ \hline
&se&\cellcolor{white}.0{\tiny 25}&0&\cellcolor{white}.0{\tiny 25}&0&\cellcolor{white}.0{\tiny 25}&\cellcolor{white}.0{\tiny 5}&\cellcolor{white}.0{\tiny 5}&\cellcolor{gray0.55}.5{\tiny 5}&\cellcolor{gray0.5}.6{\tiny 25}&\cellcolor{gray0.5}.6{\tiny 5}&\cellcolor{gray0.45}.7{\tiny 25}&\cellcolor{gray0.5}.6{\tiny 5}&\cellcolor{gray0.55}.5&\cellcolor{gray0.5}.6{\tiny 75}&\cellcolor{gray}1&\cellcolor{gray0.8}.2{\tiny 5}&\cellcolor{gray0.8}.2{\tiny 75}&\cellcolor{gray0.9}.1{\tiny 75}&\cellcolor{gray0.9}.1&\cellcolor{gray0.9}.1{\tiny 5}&\cellcolor{gray0.9}.1{\tiny 25}&\cellcolor{white}.0{\tiny 75}&\cellcolor{white}.0{\tiny 5}&\cellcolor{gray0.8}.2&\cellcolor{gray0.9}.1{\tiny 25}&\cellcolor{gray0.9}.1{\tiny 75}&\cellcolor{gray0.9}.1{\tiny 5}&\cellcolor{white}.0{\tiny 75}&\cellcolor{gray0.7}.3&\cellcolor{white}.0{\tiny 5} \\ \hline
&doubler&0&\cellcolor{white}.0{\tiny 75}&\cellcolor{white}.0{\tiny 5}&\cellcolor{white}.0{\tiny 25}&\cellcolor{white}.0{\tiny 25}&\cellcolor{white}.0{\tiny 25}&\cellcolor{white}.0{\tiny 5}&\cellcolor{gray0.8}.2&\cellcolor{white}.0{\tiny 75}&\cellcolor{white}.0{\tiny 25}&\cellcolor{gray0.9}.1{\tiny 75}&\cellcolor{white}.0{\tiny 75}&\cellcolor{gray0.8}.2{\tiny 25}&\cellcolor{gray0.6}.4{\tiny 5}&\cellcolor{gray0.8}.2{\tiny 5}&\cellcolor{gray}1&\cellcolor{gray0.35}.9{\tiny 25}&\cellcolor{gray0.55}.5{\tiny 75}&\cellcolor{gray0.7}.3{\tiny 25}&\cellcolor{gray0.5}.6{\tiny 5}&\cellcolor{gray0.7}.3{\tiny 75}&\cellcolor{gray0.9}.1&\cellcolor{white}.0{\tiny 5}&\cellcolor{gray0.8}.2{\tiny 75}&\cellcolor{white}.0{\tiny 5}&\cellcolor{white}.0{\tiny 75}&\cellcolor{white}.0{\tiny 5}&\cellcolor{white}.0{\tiny 75}&\cellcolor{white}.0{\tiny 25}&\cellcolor{gray0.9}.1 \\ \hline
&double&0&\cellcolor{gray0.9}.1&\cellcolor{white}.0{\tiny 5}&\cellcolor{white}.0{\tiny 25}&\cellcolor{white}.0{\tiny 25}&\cellcolor{white}.0{\tiny 25}&\cellcolor{white}.0{\tiny 5}&\cellcolor{gray0.8}.2&\cellcolor{gray0.9}.1{\tiny 75}&\cellcolor{white}.0{\tiny 5}&\cellcolor{gray0.9}.1{\tiny 75}&\cellcolor{gray0.9}.1&\cellcolor{gray0.7}.3{\tiny 5}&\cellcolor{gray0.6}.4{\tiny 75}&\cellcolor{gray0.8}.2{\tiny 75}&\cellcolor{gray0.35}.9{\tiny 25}&\cellcolor{gray}1&\cellcolor{gray0.55}.5&\cellcolor{gray0.8}.2{\tiny 75}&\cellcolor{gray0.5}.6{\tiny 25}&\cellcolor{gray0.7}.3{\tiny 25}&\cellcolor{gray0.9}.1&\cellcolor{white}.0{\tiny 5}&\cellcolor{gray0.8}.2{\tiny 75}&\cellcolor{white}.0{\tiny 75}&\cellcolor{white}.0{\tiny 75}&\cellcolor{white}.0{\tiny 5}&\cellcolor{white}.0{\tiny 25}&\cellcolor{white}.0{\tiny 5}&\cellcolor{white}.0{\tiny 75} \\ \hline
&partie&\cellcolor{white}.0{\tiny 5}&\cellcolor{gray0.9}.1{\tiny 5}&\cellcolor{gray0.9}.1{\tiny 25}&\cellcolor{white}.0{\tiny 25}&\cellcolor{gray0.9}.1{\tiny 25}&\cellcolor{white}.0{\tiny 5}&\cellcolor{white}.0{\tiny 75}&\cellcolor{gray0.9}.1{\tiny 25}&\cellcolor{white}.0{\tiny 25}&0&\cellcolor{gray0.9}.1&\cellcolor{white}.0{\tiny 5}&\cellcolor{gray0.8}.2&\cellcolor{gray0.7}.3&\cellcolor{gray0.9}.1{\tiny 75}&\cellcolor{gray0.55}.5{\tiny 75}&\cellcolor{gray0.55}.5&\cellcolor{gray}1&\cellcolor{gray0.55}.5{\tiny 5}&\cellcolor{gray0.5}.6{\tiny 5}&\cellcolor{gray0.45}.7{\tiny 75}&\cellcolor{gray0.8}.2{\tiny 75}&\cellcolor{gray0.8}.2{\tiny 25}&\cellcolor{gray0.7}.3{\tiny 75}&\cellcolor{gray0.8}.2{\tiny 25}&\cellcolor{gray0.9}.1{\tiny 5}&\cellcolor{gray0.9}.1{\tiny 5}&\cellcolor{gray0.9}.1{\tiny 25}&\cellcolor{white}.0{\tiny 5}&\cellcolor{gray0.9}.1 \\ \hline
&science&0&\cellcolor{white}.0{\tiny 25}&\cellcolor{white}.0{\tiny 5}&0&\cellcolor{white}.0{\tiny 25}&0&\cellcolor{white}.0{\tiny 25}&\cellcolor{white}.0{\tiny 75}&\cellcolor{white}.0{\tiny 25}&\cellcolor{white}.0{\tiny 25}&\cellcolor{gray0.9}.1{\tiny 25}&\cellcolor{gray0.9}.1{\tiny 25}&\cellcolor{gray0.8}.2&\cellcolor{white}.0{\tiny 75}&\cellcolor{gray0.9}.1&\cellcolor{gray0.7}.3{\tiny 25}&\cellcolor{gray0.8}.2{\tiny 75}&\cellcolor{gray0.55}.5{\tiny 5}&\cellcolor{gray}1&\cellcolor{gray0.6}.4{\tiny 75}&\cellcolor{gray0.5}.6{\tiny 25}&\cellcolor{gray0.7}.3{\tiny 5}&\cellcolor{gray0.8}.2{\tiny 5}&\cellcolor{gray0.9}.1{\tiny 5}&\cellcolor{gray0.9}.1&\cellcolor{white}.0{\tiny 75}&\cellcolor{white}.0{\tiny 5}&\cellcolor{white}.0{\tiny 75}&\cellcolor{white}.0{\tiny 5}&\cellcolor{gray0.9}.1 \\ \hline
&compter&0&\cellcolor{gray0.9}.1{\tiny 75}&\cellcolor{white}.0{\tiny 75}&0&\cellcolor{white}.0{\tiny 5}&0&\cellcolor{gray0.9}.1&\cellcolor{gray0.9}.1{\tiny 25}&\cellcolor{white}.0{\tiny 25}&0&\cellcolor{white}.0{\tiny 5}&\cellcolor{white}.0{\tiny 25}&\cellcolor{gray0.9}.1{\tiny 75}&\cellcolor{gray0.7}.3{\tiny 75}&\cellcolor{gray0.9}.1{\tiny 5}&\cellcolor{gray0.5}.6{\tiny 5}&\cellcolor{gray0.5}.6{\tiny 25}&\cellcolor{gray0.5}.6{\tiny 5}&\cellcolor{gray0.6}.4{\tiny 75}&\cellcolor{gray}1&\cellcolor{gray0.6}.4{\tiny 25}&\cellcolor{white}.0{\tiny 75}&\cellcolor{gray0.9}.1&\cellcolor{gray0.6}.4{\tiny 25}&\cellcolor{gray0.9}.1&\cellcolor{gray0.9}.1&\cellcolor{gray0.9}.1&\cellcolor{white}.0{\tiny 25}&\cellcolor{gray0.9}.1&\cellcolor{white}.0{\tiny 5} \\ \hline
&ensemble&\cellcolor{white}.0{\tiny 75}&\cellcolor{gray0.9}.1&\cellcolor{gray0.9}.1{\tiny 5}&\cellcolor{white}.0{\tiny 5}&\cellcolor{gray0.9}.1{\tiny 5}&\cellcolor{white}.0{\tiny 75}&0&\cellcolor{gray0.9}.1{\tiny 25}&0&0&\cellcolor{white}.0{\tiny 75}&\cellcolor{white}.0{\tiny 5}&\cellcolor{gray0.9}.1{\tiny 75}&\cellcolor{gray0.9}.1{\tiny 75}&\cellcolor{gray0.9}.1{\tiny 25}&\cellcolor{gray0.7}.3{\tiny 75}&\cellcolor{gray0.7}.3{\tiny 25}&\cellcolor{gray0.45}.7{\tiny 75}&\cellcolor{gray0.5}.6{\tiny 25}&\cellcolor{gray0.6}.4{\tiny 25}&\cellcolor{gray}1&\cellcolor{gray0.6}.4{\tiny 25}&\cellcolor{gray0.8}.2{\tiny 5}&\cellcolor{gray0.9}.1{\tiny 5}&\cellcolor{gray0.8}.2{\tiny 75}&\cellcolor{white}.0{\tiny 75}&\cellcolor{gray0.9}.1&\cellcolor{white}.0{\tiny 5}&\cellcolor{white}.0{\tiny 25}&\cellcolor{gray0.9}.1 \\ \hline
&dire&0&\cellcolor{white}.0{\tiny 25}&\cellcolor{gray0.9}.1{\tiny 25}&0&\cellcolor{white}.0{\tiny 75}&\cellcolor{white}.0{\tiny 25}&\cellcolor{white}.0{\tiny 25}&\cellcolor{gray0.9}.1{\tiny 25}&\cellcolor{white}.0{\tiny 5}&\cellcolor{white}.0{\tiny 5}&\cellcolor{gray0.8}.2&\cellcolor{gray0.8}.2{\tiny 25}&\cellcolor{gray0.9}.1&\cellcolor{white}.0{\tiny 25}&\cellcolor{white}.0{\tiny 75}&\cellcolor{gray0.9}.1&\cellcolor{gray0.9}.1&\cellcolor{gray0.8}.2{\tiny 75}&\cellcolor{gray0.7}.3{\tiny 5}&\cellcolor{white}.0{\tiny 75}&\cellcolor{gray0.6}.4{\tiny 25}&\cellcolor{gray}1&\cellcolor{gray0.7}.3&\cellcolor{white}.0{\tiny 25}&\cellcolor{white}.0{\tiny 5}&0&0&\cellcolor{white}.0{\tiny 25}&\cellcolor{white}.0{\tiny 25}&\cellcolor{white}.0{\tiny 25} \\ \hline
&regle&\cellcolor{gray0.9}.1{\tiny 25}&\cellcolor{white}.0{\tiny 75}&\cellcolor{gray0.9}.1&\cellcolor{white}.0{\tiny 5}&\cellcolor{gray0.9}.1&\cellcolor{white}.0{\tiny 5}&\cellcolor{gray0.9}.1{\tiny 75}&\cellcolor{gray0.9}.1&\cellcolor{white}.0{\tiny 25}&0&0&\cellcolor{white}.0{\tiny 5}&\cellcolor{white}.0{\tiny 5}&\cellcolor{white}.0{\tiny 25}&\cellcolor{white}.0{\tiny 5}&\cellcolor{white}.0{\tiny 5}&\cellcolor{white}.0{\tiny 5}&\cellcolor{gray0.8}.2{\tiny 25}&\cellcolor{gray0.8}.2{\tiny 5}&\cellcolor{gray0.9}.1&\cellcolor{gray0.8}.2{\tiny 5}&\cellcolor{gray0.7}.3&\cellcolor{gray}1&\cellcolor{gray0.8}.2{\tiny 25}&\cellcolor{gray0.8}.2{\tiny 75}&\cellcolor{gray0.8}.2{\tiny 75}&\cellcolor{gray0.7}.3{\tiny 25}&\cellcolor{gray0.7}.3{\tiny 25}&\cellcolor{gray0.9}.1&\cellcolor{gray0.9}.1{\tiny 5} \\ \hline
&chose&\cellcolor{gray0.9}.1&\cellcolor{gray0.9}.1{\tiny 75}&\cellcolor{gray0.9}.1&\cellcolor{white}.0{\tiny 75}&\cellcolor{gray0.9}.1{\tiny 75}&\cellcolor{white}.0{\tiny 75}&\cellcolor{white}.0{\tiny 75}&\cellcolor{white}.0{\tiny 5}&\cellcolor{gray0.8}.2{\tiny 25}&\cellcolor{white}.0{\tiny 25}&\cellcolor{gray0.9}.1&\cellcolor{white}.0{\tiny 25}&\cellcolor{gray0.8}.2&\cellcolor{gray0.9}.1{\tiny 5}&\cellcolor{gray0.8}.2&\cellcolor{gray0.8}.2{\tiny 75}&\cellcolor{gray0.8}.2{\tiny 75}&\cellcolor{gray0.7}.3{\tiny 75}&\cellcolor{gray0.9}.1{\tiny 5}&\cellcolor{gray0.6}.4{\tiny 25}&\cellcolor{gray0.9}.1{\tiny 5}&\cellcolor{white}.0{\tiny 25}&\cellcolor{gray0.8}.2{\tiny 25}&\cellcolor{gray}1&\cellcolor{gray0.7}.3{\tiny 75}&\cellcolor{gray0.6}.4{\tiny 75}&\cellcolor{gray0.6}.4{\tiny 25}&\cellcolor{gray0.55}.5&\cellcolor{gray0.9}.1{\tiny 5}&\cellcolor{gray0.8}.2{\tiny 5} \\ \hline
&reduire&\cellcolor{gray0.8}.2&\cellcolor{gray0.9}.1{\tiny 5}&\cellcolor{gray0.8}.2&\cellcolor{gray0.9}.1{\tiny 5}&\cellcolor{gray0.8}.2{\tiny 75}&\cellcolor{gray0.7}.3{\tiny 25}&\cellcolor{white}.0{\tiny 75}&\cellcolor{white}.0{\tiny 5}&\cellcolor{gray0.8}.2{\tiny 25}&\cellcolor{white}.0{\tiny 5}&\cellcolor{gray0.9}.1&0&\cellcolor{gray0.7}.3{\tiny 25}&\cellcolor{white}.0{\tiny 75}&\cellcolor{gray0.9}.1{\tiny 25}&\cellcolor{white}.0{\tiny 5}&\cellcolor{white}.0{\tiny 75}&\cellcolor{gray0.8}.2{\tiny 25}&\cellcolor{gray0.9}.1&\cellcolor{gray0.9}.1&\cellcolor{gray0.8}.2{\tiny 75}&\cellcolor{white}.0{\tiny 5}&\cellcolor{gray0.8}.2{\tiny 75}&\cellcolor{gray0.7}.3{\tiny 75}&\cellcolor{gray}1&\cellcolor{gray0.8}.2{\tiny 5}&\cellcolor{gray0.8}.2{\tiny 75}&\cellcolor{gray0.9}.1&\cellcolor{gray0.8}.2{\tiny 75}&\cellcolor{white}.0{\tiny 75} \\ \hline
&seulement&\cellcolor{white}.0{\tiny 25}&\cellcolor{white}.0{\tiny 25}&0&0&\cellcolor{white}.0{\tiny 25}&0&\cellcolor{gray0.9}.1{\tiny 25}&\cellcolor{white}.0{\tiny 25}&\cellcolor{gray0.9}.1{\tiny 25}&\cellcolor{white}.0{\tiny 25}&\cellcolor{white}.0{\tiny 5}&\cellcolor{white}.0{\tiny 5}&\cellcolor{white}.0{\tiny 5}&\cellcolor{gray0.9}.1{\tiny 75}&\cellcolor{gray0.9}.1{\tiny 75}&\cellcolor{white}.0{\tiny 75}&\cellcolor{white}.0{\tiny 75}&\cellcolor{gray0.9}.1{\tiny 5}&\cellcolor{white}.0{\tiny 75}&\cellcolor{gray0.9}.1&\cellcolor{white}.0{\tiny 75}&0&\cellcolor{gray0.8}.2{\tiny 75}&\cellcolor{gray0.6}.4{\tiny 75}&\cellcolor{gray0.8}.2{\tiny 5}&\cellcolor{gray}1&\cellcolor{gray0.45}.7{\tiny 75}&\cellcolor{gray0.5}.6{\tiny 75}&\cellcolor{gray0.7}.3{\tiny 5}&\cellcolor{white}.0{\tiny 5} \\ \hline
&pratique&0&\cellcolor{white}.0{\tiny 5}&\cellcolor{white}.0{\tiny 25}&0&0&0&\cellcolor{gray0.9}.1&\cellcolor{gray0.9}.1&\cellcolor{gray0.9}.1&\cellcolor{white}.0{\tiny 25}&\cellcolor{white}.0{\tiny 5}&\cellcolor{white}.0{\tiny 5}&\cellcolor{white}.0{\tiny 75}&\cellcolor{gray0.9}.1{\tiny 25}&\cellcolor{gray0.9}.1{\tiny 5}&\cellcolor{white}.0{\tiny 5}&\cellcolor{white}.0{\tiny 5}&\cellcolor{gray0.9}.1{\tiny 5}&\cellcolor{white}.0{\tiny 5}&\cellcolor{gray0.9}.1&\cellcolor{gray0.9}.1&0&\cellcolor{gray0.7}.3{\tiny 25}&\cellcolor{gray0.6}.4{\tiny 25}&\cellcolor{gray0.8}.2{\tiny 75}&\cellcolor{gray0.45}.7{\tiny 75}&\cellcolor{gray}1&\cellcolor{gray0.5}.6{\tiny 25}&\cellcolor{gray0.7}.3{\tiny 75}&\cellcolor{white}.0{\tiny 5} \\ \hline
&multiplication&\cellcolor{white}.0{\tiny 75}&0&0&0&\cellcolor{white}.0{\tiny 5}&\cellcolor{white}.0{\tiny 5}&\cellcolor{white}.0{\tiny 5}&0&\cellcolor{gray0.9}.1&\cellcolor{white}.0{\tiny 25}&0&\cellcolor{white}.0{\tiny 25}&\cellcolor{white}.0{\tiny 25}&\cellcolor{white}.0{\tiny 25}&\cellcolor{white}.0{\tiny 75}&\cellcolor{white}.0{\tiny 75}&\cellcolor{white}.0{\tiny 25}&\cellcolor{gray0.9}.1{\tiny 25}&\cellcolor{white}.0{\tiny 75}&\cellcolor{white}.0{\tiny 25}&\cellcolor{white}.0{\tiny 5}&\cellcolor{white}.0{\tiny 25}&\cellcolor{gray0.7}.3{\tiny 25}&\cellcolor{gray0.55}.5&\cellcolor{gray0.9}.1&\cellcolor{gray0.5}.6{\tiny 75}&\cellcolor{gray0.5}.6{\tiny 25}&\cellcolor{gray}1&\cellcolor{white}.0{\tiny 5}&\cellcolor{gray0.7}.3{\tiny 25} \\ \hline
&exemple&\cellcolor{white}.0{\tiny 5}&\cellcolor{white}.0{\tiny 5}&\cellcolor{white}.0{\tiny 5}&\cellcolor{white}.0{\tiny 25}&\cellcolor{white}.0{\tiny 5}&\cellcolor{white}.0{\tiny 75}&\cellcolor{gray0.9}.1{\tiny 25}&\cellcolor{gray0.7}.3{\tiny 25}&\cellcolor{gray0.9}.1{\tiny 5}&\cellcolor{gray0.9}.1{\tiny 5}&\cellcolor{white}.0{\tiny 75}&\cellcolor{white}.0{\tiny 75}&\cellcolor{gray0.9}.1&\cellcolor{gray0.7}.3&\cellcolor{gray0.7}.3&\cellcolor{white}.0{\tiny 25}&\cellcolor{white}.0{\tiny 5}&\cellcolor{white}.0{\tiny 5}&\cellcolor{white}.0{\tiny 5}&\cellcolor{gray0.9}.1&\cellcolor{white}.0{\tiny 25}&\cellcolor{white}.0{\tiny 25}&\cellcolor{gray0.9}.1&\cellcolor{gray0.9}.1{\tiny 5}&\cellcolor{gray0.8}.2{\tiny 75}&\cellcolor{gray0.7}.3{\tiny 5}&\cellcolor{gray0.7}.3{\tiny 75}&\cellcolor{white}.0{\tiny 5}&\cellcolor{gray}1&0 \\ \hline
&nombrer&\cellcolor{white}.0{\tiny 25}&\cellcolor{white}.0{\tiny 25}&\cellcolor{white}.0{\tiny 25}&\cellcolor{white}.0{\tiny 75}&\cellcolor{white}.0{\tiny 25}&\cellcolor{white}.0{\tiny 25}&0&0&\cellcolor{white}.0{\tiny 75}&0&0&\cellcolor{white}.0{\tiny 25}&\cellcolor{white}.0{\tiny 75}&0&\cellcolor{white}.0{\tiny 5}&\cellcolor{gray0.9}.1&\cellcolor{white}.0{\tiny 75}&\cellcolor{gray0.9}.1&\cellcolor{gray0.9}.1&\cellcolor{white}.0{\tiny 5}&\cellcolor{gray0.9}.1&\cellcolor{white}.0{\tiny 25}&\cellcolor{gray0.9}.1{\tiny 5}&\cellcolor{gray0.8}.2{\tiny 5}&\cellcolor{white}.0{\tiny 75}&\cellcolor{white}.0{\tiny 5}&\cellcolor{white}.0{\tiny 5}&\cellcolor{gray0.7}.3{\tiny 25}&0&\cellcolor{gray}1 \\ \hline
\end{tabular}
}
\caption{Same as Table~\ref{bertinfickle1}, reorganized and with shades proportional to value.}
\label{bertin4}
\end{table}

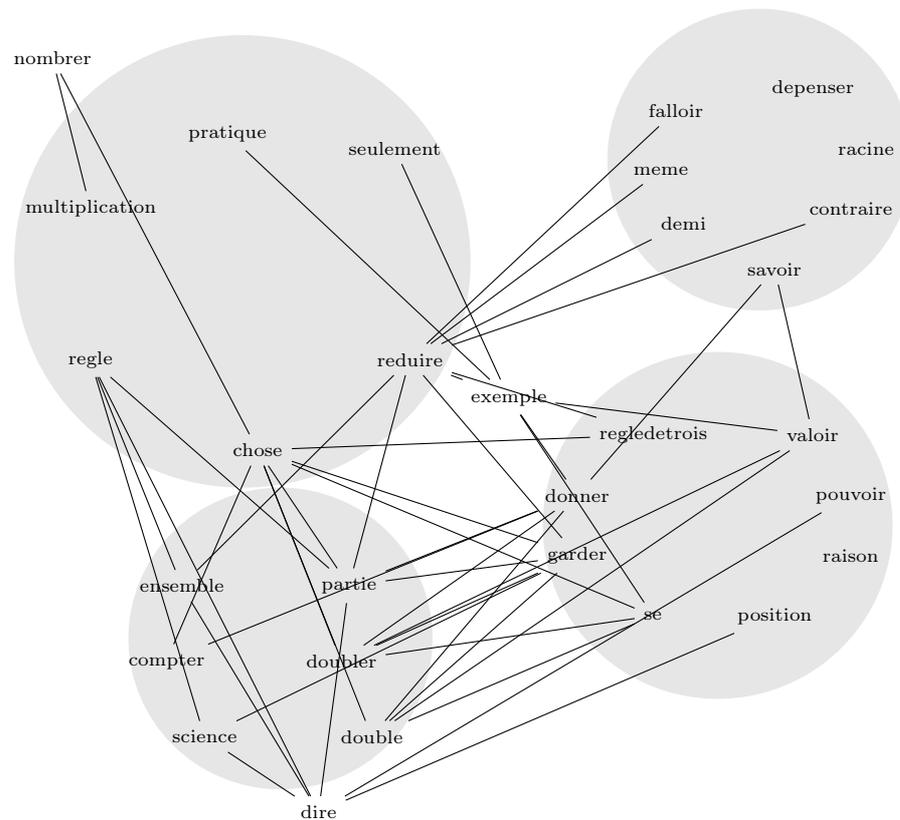
\begin{figure}[h!]
\centering
\begin{tikzpicture}[scale=1]

\fill[grisclair] (demi)+(1,0.85) circle(2);
\fill[grisclair] (se)+(0.85,1.2) circle(2.3);
\fill[grisclair] (regle)+(2.5,1.3) circle(3);
\fill[grisclair] (science)+(1,0.3) circle(2);

\begin{scope}[every node/.style={font={\scriptsize}}]

\node(multiplication) at (1,8){multiplication};
\node(seulement) at (5,8.8){seulement};
\node(pratique) at (2.8,9){pratique};
\node(regle) at (1,6){regle};
\node(reduire) at (5.2,6){reduire};
\node(chose) at (3.2,4.8){chose};
\node(depenser) at (10.5,9.6){depenser};
\node(falloir) at (8.7,9.3){falloir};
\node(meme) at (8.5,8.5){meme};
\node(racine) at (11.2,8.8){racine};
\node(contraire) at (11,8){contraire};
\node(demi) at (8.8,7.8){demi};
\node(savoir) at (10,7.2){savoir};
\node(compter) at (2,2){compter};
\node(ensemble) at (2.2,3){ensemble};
\node(partie) at (4.4,3){partie};
\node(double) at (4.7,1){double};
\node(doubler) at (4.3,2){doubler};
\node(science) at (2.5,1){science};
\node(garder) at (7.4,3.4){garder};
\node(donner) at (7.4,4.2){donner};
\node(pouvoir) at (11,4.2){pouvoir};
\node(raison) at (11,3.4){raison};
\node(regledetrois) at (8.4,5){regledetrois};
\node(se) at (8.4,2.6){se};
\node(valoir) at (10.5,5){valoir};
\node(position) at (10,2.6){position};
\node(dire) at (4,0){dire};
\node(exemple) at (6.5,5.5){exemple};
\node(nombrer) at (0.5,10){nombrer};
\end{scope}

\draw(compter) -- (chose);
\draw(donner) -- (compter);
\draw(double) -- (chose);
\draw(double) -- (donner);
\draw(doubler) -- (chose);
\draw(doubler) -- (donner);
\draw(ensemble) -- (dire);
\draw(exemple) -- (donner);
\draw(garder) -- (chose);
\draw(garder) -- (double);
\draw(garder) -- (doubler);
\draw(nombrer) -- (chose);
\draw(nombrer) -- (multiplication);
\draw(partie) -- (chose);
\draw(partie) -- (dire);
\draw(partie) -- (donner);
%\draw(partie) -- (figurer);
\draw(partie) -- (garder);
\draw(position) -- (dire);
%\draw(position) -- (figurer);
\draw(pouvoir) -- (dire);
%\draw(pouvoir) -- (figurer);
\draw(pratique) -- (exemple);
\draw(reduire) -- (contraire);
\draw(reduire) -- (demi);
\draw(reduire) -- (ensemble);
\draw(reduire) -- (exemple);
\draw(reduire) -- (falloir);
\draw(reduire) -- (garder);
\draw(reduire) -- (meme);
\draw(reduire) -- (partie);
\draw(regle) -- (dire);
\draw(regle) -- (ensemble);
%\draw(regle) -- (figurer);
\draw(regle) -- (partie);
\draw(regledetrois) -- (chose);
\draw(regledetrois) -- (reduire);
\draw(savoir) -- (donner);
\draw(science) -- (dire);
%\draw(science) -- (figurer);
\draw(science) -- (garder);
\draw(science) -- (regle);
\draw(se) -- (chose);
\draw(se) -- (double);
\draw(se) -- (doubler);
\draw(se) -- (exemple);
\draw(seulement) -- (exemple);
\draw(valoir) -- (double);
\draw(valoir) -- (doubler);
\draw(valoir) -- (exemple);
\draw(valoir) -- (savoir);
\end{tikzpicture}
\caption{Glutton decomposition in quasi-cliques of maximum size.}
\label{fig:graph2}
\end{figure}

\subsection{Neighborhood graphs}

Bertin matrix (see Table~\ref{bertin4}) shows some remarkable clustering among fickle words. The question is now to produce some interpretation of this clustering.\\

First, we can observe that the Bertin matrix displays four groups along its diagonal, from the more connected on the top left, to the less connected on the bottom right.\\

The first list (\textit{contraire}, \textit{depenser}, \textit{falloir}, \textit{racine}, \textit{meme}, \textit{demi}, \textit{savoir}, see Figure~\ref{fig:ficklest} for a translation) is a collection of rather heterogeneous words. There are words frequently used such as \textit{demi} and others bearing a strong polysemy such as \textit{falloir}.\\
 
A possible explanation may be that these groups of words form phrases in the corpus (that is usually called \textit{textual co-occurence): the words are spaced only one or two words  from each other.} and that these associations are reflected in the table. In that case, fickle words can be used as a tool to extract \textit{topoi}. Thus, for example, \textit{savoir} "to know" and \textit{contraire} "contrary" are often used in phrases such as \textit{savoir par son contraire} "to know smth through its contrary".\\

We can also think that these clustering properties reveal more distant co-occurrences, that means words appearing in the same sentence or paragraph, but not necessarily the same \textit{topos}. For instance, \textit{falloir} "to have to do" has a lot of such co-occurrences with fickle words like \textit{reduire} "to reduce", \textit{racine} "root" and \textit{savoir} "to know". In this configuration, we can in fact suppose that all these words are shared by the same sentences.\\

On the contrary, \textit{demi}, that is a part of the same well-connected group (according to the clustering), does not have any specific co-occurrences in the texts with any word in this group. That is especially interesting, since it reveals the existence of connections that could not have been deduced from a simple study of co-occurrence with classical tools.\\

The Figure~\ref{fig:graph2} opens another perspective. Indeed, it shows the words that make the link between clusters. These fickle words have a lot of different affinities. We can see, for example, that the positions of \textit{reduire} "to reduce" and \textit{exemple} "example" are not very surprising, because these words are used a lot, in every text, in sentences associating them with various other fickle words, such as "in all the examples preceding the problems", or "the problem of reducing or converting the monetary values".

These questions are not solved yet, and the answer cannot be sure without an enlargement of the corpus. Indeed, we would like to test this hypothesis by using the process described here on a larger part of the corpus.

\section*{Conclusion}

In this work, we have shown how to use the Kohonen maps as a complement of Factorial Correspondence  Analysis  methods (FCA)classically used in lexicometry,

\begin{itemize}
\item to improve the information provided by the different projections of the FCA,

\item to make the Kohonen maps more robust with respect to the randomness of the SOM  algorithm, by distinguishing stable neighbor pairs from fickle pairs,

\item to build graphs of connections between fickle words which are difficult to analyze by both FCA and Kohonen map alone.
\end{itemize}

We think that it will be interesting to use this methodology on a large variety of corpus, such as political speeches, chivalric culture~\cite{Der2010} texts and scientific articles. 

%\section*{References}

\bibliography{wsom2012_bis}

\end{document}